\documentclass[reqno,11pt,twoside,english]{article}
\usepackage{amsmath, amsfonts, amssymb, amsthm, amscd}
\usepackage[T1]{fontenc}
\usepackage[english]{babel}

\DeclareMathSymbol{\leqslant}{\mathalpha}{AMSa}{"36} 
\DeclareMathSymbol{\geqslant}{\mathalpha}{AMSa}{"3E} 
\renewcommand{\leq}{\;\leqslant\;}                   
\renewcommand{\geq}{\;\geqslant\;}                   

\newtheorem{theorem}{{\bf Theorem}}[section]
\newtheorem{lemma}{{\bf Lemma}}[section]
\newtheorem{proposition}{{\bf Proposition}}[section]
\newtheorem{definition}{{\bf Definition}}[section]

\newtheorem{remark}{{\bf Remark}}[section]

\newtheorem{assumption}{{\bf Assumption}}[section]

\newcommand{\be}{\begin{equation}}
\newcommand{\ee}{\end{equation}}
\newcommand{\beq}{\begin{eqnarray}}
\newcommand{\eeq}{\end{eqnarray}}

\newcommand{\ced}{\end{proof}}

\def\first {\vtop{\baselineskip=11pt
\hbox to 125truept{\hss Anis MATOUSSI\hss}\vskip2truept \hbox
to 125truept{\small\hss Equipe ``Statistique et Processus''\hss}
\hbox to 125truept{\small\hss Universit\'e du Maine\hss} \hbox to
125truept{\small\hss Avenue Olivier Messiaen\hss} \hbox to
125truept{\small\hss 72085 LE MANS Cedex 9 \hss} \hbox to
125truept{\small\hss FRANCE\hss}
\hbox to 125truept{\small\hss
anis.matoussi@univ-lemans.fr\hss}}}
\def\second {\vtop{\baselineskip=11pt
\hbox to 125truept{\hss Mingyu XU\hss }\vskip2truept \hbox to
125truept{\small\hss Institute of Applied Mathematics \hss} \hbox to
125truept{\small\hss Academy of Mathematics and Systems Science\hss} \hbox to
125truept{\small\hss CAS, Beijing, 100190  China \hss} \hbox to
125truept{\small\hss xumy@amss.ac.cn\hss }}}
\begin{document}
\author{\first \and \second }
\title{Sobolev solution for semilinear PDE with obstacle under monotonicity condition}
\date{\it{\footnotesize first revision version  March 20, 2006, final accepted version May 26, 2008}}
 \maketitle

 \abstract{We prove the
existence and uniqueness of Sobolev solution of a semilinear PDE's
and  PDE's with obstacle under monotonicity condition. Moreover we
give the probabilistic interpretation of the solutions in term of
Backward SDE and reflected Backward SDE respectively.}
\\[0.5cm]
\textbf{Key words:} Backward stochastic differential equation,
Reflected backward stochastic differential equation, monotonicity
condition, Stochastic flow, partial differential equation with
obstacle.
\\[0.5cm]
\textbf{AMS Classification:} 35D05, 60H10, 60H30B
\\[0.5cm]

\section{Introduction}
Our approach is  based on Backward Stochastic Differential Equations
(in short BSDE's) which were first introduced by Bismut
\cite{Bismut} in 1973 as equation for the adjoint process in the
stochastic version of Pontryagin maximum principle. Pardoux and Peng
\cite{PP90} generalized the notion in 1990 and were the first to
consider general BSDE's and to solve the question of existence and
uniqueness in the non-linear case. Since then BSDE's have been
widely used in stochastic control and especially in mathematical
finance, as any pricing problem by replication can be written in
terms of linear BSDEs, or non-linear BSDEs when portfolios
constraints are taken into account as in El Karoui,  Peng and Quenez
\cite{EPQ}.

The main motivation  to introduce the non-linear BSDE's was to give
a probabilistic interpretation (Feynman-Kac's formula) for the
solutions of semilinear parabolic PDE's. This result was first
obtained by Peng
 in \cite{P91}, see also Pardoux and Peng
\cite{PP92} by considering the viscosity and classical  solutions of
such PDE's. Later, Barles and Lesigne
 \cite{BL} studied the relation between BSDE's and solutions of
semi-linear PDE's in Soblev spaces. More recently Bally and Matoussi
\cite{BM2001} studied semilinear stochastic PDEs and backward doubly
SDE in Sobolev space and  their probabilistic method is based on
stochastic flow.

The reflected BSDE's was introduced by the five authors El Karoui,
Kapoudjian, Pardoux,  Peng and Quenez in \cite{EKPPQ}, the setting
of those equations is the following: let us consider moreover an
adapted stochastic process $L:=(L_t)_{t\leq T}$ which stands for a
barrier. A solution for the reflected BSDE associated with
$(\xi,g,L)$ is a triple of adapted stochastic processes
$(Y_t,Z_t,K_t)_{t\leq T}$ such that
\begin{equation*}
\left\lbrace
\begin{split}
&  Y_t=\xi+ \displaystyle \int_{t}^{T}g(s,\omega,
Y_s,Z_s)ds+K_T-K_t-\int_{t}^{T}Z_sdB_s, \,\,\forall \,  t\,  \in [0,  T],\\
& \displaystyle Y_t\geq L_t \; \mbox { and } \;
\int_0^T(Y_t-L_t)dK_t=0.
\end{split}
\right.
\end{equation*} The process $K$
is continuous, increasing and its role is to push upward $Y$ in
order to keep it above the barrier $L$. The requirement
$\int_0^T(Y_t-L_t)dK_t=0$ means that the action of $K$ is made with
a minimal energy.

The development of reflected BSDE's (see for example \cite{EKPPQ},
\cite{HLM}, \cite{HO}) has been especially motivated by pricing
American contingent claim by replication, especially in constrained
markets. Actually it has been shown by El Karoui, Pardoux and Quenez
\cite{EKPQ} that the price of an American contingent claim
$(S_t)_{t\leq T}$ whose strike is $\gamma$ in a standard complete
financial market is $Y_0$ where $(Y_t,\pi_t,K_t)_{t\leq T}$ is the
solution of the following reflected BSDE
$$ \left\{
\begin{array}{l}
-dY_t=b(t,Y_t,\pi_t)dt+dK_t-\pi_tdW_t,\quad Y_T=(S_T-\gamma)^+,\\
 Y_t\geq (S_t-\gamma)^+ \quad \mbox{ and }\quad
\int_0^T(Y_t-(S_t-\gamma)^+)dK_t=0
\end{array}
\right. $$ for an appropriate choice of the function $b$. The
process $\pi$ allows to construct a replication strategy and $K$ is
a consumption process that could have the buyer of the option. In a
standard financial market the function $b(t,\omega,y,z)=
r_ty+z\theta_t$ where $\theta_t$ is the risk premium and $r_t$ the
spot rate to invest or borrow. Now when the market is constrained
$i.e.$ the interest rates are not the same whether we borrow or
invest money then the function $b(t,\omega,y,z)=
r_ty+z\theta_t-(R_t-r_t)(y-(z.\sigma_t^{-1}.\bf{1}))^-$ where $R_t$
(resp. $r_t$) is the spot rate to borrow (resp. invest) and $\sigma$
the volatility.

Partial Differential Equations with obstacles and their connections
with optimal control problems have been studied by Bensoussan and
Lions \cite{BenL}.   They study such equations in the point of view
of variational inequalities. In a recent paper, Bally, Caballero, El
Karoui and Fernandez \cite{BCEF} studied the the following
semilinear PDE with obstacle
\[
(\partial _{t}+\mathcal{L})u+f(t,x,u,\sigma ^{*}\nabla u)+\nu
=0,\;u\geq h,\;u_{T}=g,
\]
where $h$ is the obstacle. The solution of such equation is a pair
$(u,\nu )$ where $u$ is a function in
$\mathbf{L}^{2}([0,T],\mathcal{H})$ and $\nu $ is a positive
measure concentrated on the set $\{u=h\}$. The authors proved the
uniqueness and existence for the solution to this PDE when the
coefficient $f $ is Lipschitz and linear increasing on $(y,z)$,
and gave the probabilistic
interpretation (Feynman-Kac formula) for $u$ and $\nabla u$ by the solution $%
(Y,Z)$ of the reflected BSDE (in short RBSDE). They prove also the
natural relation between Reflected BSDE's and variational
inequalities and prove uniqueness of the solution for such
variational problem  by using the relation between the increasing
process $K$ and the measure $\nu $. This is also a point of view in
this paper.

On the other hand, Pardoux \cite{Pardoux99}  studied the solution of
a BSDE with a coefficient $f(t,\omega ,y,z)$, which satisfies only
monotonicity, continuous and general increasing conditions on $y$,
and a Lipschitz
condition on $z$, i.e. for some continuous, increasing function $\varphi :%
\mathbb{R}_{+}\rightarrow \mathbb{R}_{+}$, and real numbers $\mu \in \mathbb{%
R}$, $k>0$, $\forall t \in [0,T]$, $\forall y,y^{\prime }\in
\mathbb{R}^{n}$, $\forall z,z^{\prime }\in \mathbb{R}^{n\times
d}$,
\begin{eqnarray}
\left| f(t,y,0)\right| &\leq &\left| f(t,0,0)\right| +\varphi
(\left| y\right| )\text{, a.s.};
\label{assup} \\
\left\langle y-y^{\prime },f(t,y,z)-f(t,y^{\prime
},z)\right\rangle &\leq &\mu \left| y-y^{\prime }\right|
^{2}\text{, a.s.};
\nonumber \\
\left| f(t,y,z)-f(t,y,z^{\prime })\right| &\leq &k\left|
z-z^{\prime }\right| \text{, a.s.}.  \nonumber
\end{eqnarray}
In the same paper, he also considered the PDE whose coefficient
$f$ satisfies the monotonicity condition (\ref{assup}), proved the
existence of a viscosity solution $u$ to this PDE and gave its
probabilistic interpretation via the solution of the corresponding
BSDE. More recently, Lepeltier, Matoussi and Xu \cite{LMX} proved
the existence and uniqueness of the solution for the reflected
BSDE under the monotonicity condition.

In our paper, we study the Sobolev solutions of the PDE and also the
PDE with continuous obstacle under the monotonicity condition
(\ref{assup}). Using penalization method, we prove the existence of
the solution and give the probabilistic interpretation of
the solution $u$ and $\nabla u$ (resp.$(u,\nabla u,\nu )$) by the solution $%
(Y,Z)$ of backward SDE (resp. the solution $(Y,Z,K)$ of reflected
backward SDE). Furthermore we use equivalence norm results and a
stochastic test function to pass from the solution of PDE's to the
one of BSDE's in order to get the uniqueness of the solution.

Our paper is organized as following: in section 2, we present the
basic assumptions and the definitions of the solutions for PDE and
PDE with obstacle, then in section 3, we recall some useful
results from \cite{BM2001}. We will prove the main results for PDE
and PDE with continuous barrier under monotonicity condition in
section 4 and 5 respectively. Finally, we prove an analogue result
to Proposition 2.3 in \cite{BM2001}  under the monotonicity
condition, and we also give a priori estimates for the solution of
the  reflected BSDE's.

\section{Notations and preliminaries}
\sectionmark{Notations and preliminaries}

Let $(\Omega ,\mathcal{F},P)$ be a complete probability space, and $%
B=(B_{1},B_{2},\cdots ,B_{d})^{*}$ be a $d$-dimensional Brownian motion
defined on a finite interval $[0,T]$, $0<T<+\infty $. Denote by $\{\mathcal{F%
}_{s}^{t};t\leq s\leq T\}$ the natural filtration generated by the Brownian
motion $B:$%
\[
\mathcal{F}_{s}^{t}=\sigma \{B_{s}-B_{t};t\leq r\leq s\}\cup \mathcal{F}%
_{0},
\]
where $\mathcal{F}_{0}$ contains all $P-$null sets of $\mathcal{F}$.

We will need the following spaces for studying BSDE or reflected BSDE. For
any given $n\in \mathbb{N}$:

\begin{itemize}
\item  $\mathbf{L}_{n}^{2}(\mathcal{F}_{s}^{t}):$ the set of $n$-dimensional
$\mathcal{F}_{s}^{t}$-measurable random variable $\xi $, such that $E(|\xi
|^{2})<+\infty $.

\item  $\mathbf{H}_{n\times m}^{2}(t,T):$ the set of
$\mathbb{R}^{m\times n}$-valued $\mathcal{F}_{s}^{t}$-predictable
process $\psi $ on the interval $[t,T]$, such that
$E\int_{t}^{T}\left\| \psi (s)\right\| ^{2}ds<+\infty .$

\item  $\mathbf{S}_{n}^{2}(t,T):$ the set of $n$%
-dimensional $\mathcal{F}_{s}^{t}$-progressively measurable process $\psi $
on the interval $[t,T]$, such that $E(\sup_{t\leq s\leq T}\left\| \psi
(s)\right\| ^{2})<+\infty .$

\item  $\mathbf{A}^{2}(t,T):=$\{$K$ : $\Omega \times [t,T]\rightarrow %
\mathbb{R}$, $\mathcal{F}_{s}^t$--progressively measurable
increasing RCLL processes \newline\hspace*{1.5cm} with $K_t=0$,
$E[(K_{T})^{2}]<\infty $ \}.
\end{itemize}

Finally, we shall denote by $\mathcal{P}$ the $\sigma $-algebra of
predictable sets on $[0,T]\times \Omega $. In the real--valued case, i.e., $%
n=1$, these spaces will be simply denoted by $\mathbf{L}^{2}(\mathcal{F}%
_{s}^{t})$, $\mathbf{H}^{2}(t,T)$ and $\mathbf{S}^{2}(t,T)$,
respectively.

For the sake of the Sobolev solution of the PDE, the following notations are
needed:

\begin{itemize}
\item  $C_{b}^{m}(\mathbb{R}^{d},\mathbb{R}^{n}):$ the set of $C^{m}$%
-functions $f:\mathbb{R}^{d}\rightarrow \mathbb{R}^{n}$, whose
partial derivatives of order less that or equal to $m$, are
bounded. (The functions themselves need not to be bounded)

\item  $C_{c}^{1,m}([0,T]\times \mathbb{R}^{d},\mathbb{R}^{n}):$
the set of continuous functions $f:[0,T]\times
\mathbb{R}^{d}\rightarrow \mathbb{R}^{n} $ with compact support,
whose first partial derivative with respect to $t$ and partial
derivatives of order less or equal to $m$ with respect to $x$
exist.

\item  $\rho :\mathbb{R}^{d}\rightarrow \mathbb{R}$, the weight, is a
continuous positive function which satisfies $\int_{%
\mathbb{R}^{d}}\rho (x)dx<\infty $.

\item  $\mathbf{L}^{2}(\mathbb{R}^{d},\rho (x)dx):$ the weighted $\mathbf{L}%
^{2}$-space with weight function $\rho (x)$, endowed with the norm
\[
\left\| u\right\|^2 _{\mathbf{L}^{2}(\mathbb{R}^{d},\rho )}=\int_{\mathbb{R}%
^{d}}\left| u(x)\right| ^{2}\rho (x)dx
\]
\end{itemize}

We assume:

\begin{assumption}\label{term-pde}$g(\cdot
)\in \mathbf{L}^{2}\mathbb{(R}^{d},\rho (x)dx)$.
\end{assumption}
\begin{assumption}\label{coef-pde} $f:[0,T]\times \mathbb{R}%
^{d}\times \mathbb{R}^{n}\mathbb{\times R}^{n\times d}\rightarrow \mathbb{R}%
^{n}$ is measurable in $(t,x,y,z)$ and
\[
\int_{0}^{T}\int_{\mathbb{R}^{d}}\left| f(t,x,0,0)\right| ^{2}\rho
(x)dxdt<\infty .
\]
\end{assumption}
\begin{assumption}\label{coef-pde1}$f$ satisfies
increasing and
monotonicity condition on $y$, for some continuous increasing function $%
\varphi :\mathbb{R}_{+}\rightarrow \mathbb{R}_{+}$, real numbers $k>0$, $\mu
\in \mathbb{R}$ such that $\forall (t,x,y,y^{\prime },z,z^{\prime })\in
[0,T]\times \mathbb{R}^{d}\times \mathbb{R}^{n}\times \mathbb{R}^{n}\times %
\mathbb{R}^{n\times d}\times \mathbb{R}^{n\times d}$

\begin{enumerate}
\item[(i)]  $\left| f(t,x,y,z)\right| \leq \left| f(t,x,0,z)\right| +\varphi
(\left| y\right| )$,

\item[(ii)]  $\left| f(t,x,y,z)-f(t,x,y,z^{\prime })\right| \leq k\left|
z-z^{\prime }\right| $,

\item[(iii)]  $\left\langle y-y^{\prime },f(t,x,y,z)-f(t,x,y^{\prime
},z)\right\rangle \leq \mu \left| y-y^{\prime }\right| ^{2}$,

\item[(iv)]  $y\rightarrow f(t,x,y,z)$ is continuous.
\end{enumerate}
\end{assumption}

For the PDE with obstacle, we consider that $f$ satisfies assumptions \ref
{coef-pde} and \ref{coef-pde1}, for $n=1$.\\[0.1cm]
\begin{assumption}\label{bar-pde} The obstacle function $h\in
C([0,T]\times \mathbb{R}^{d};\mathbb{R})$ satisfies the following
conditions: there exists $\kappa \in \mathbb{R}$, $\beta >0$, such
that $\forall (t,x)\in [0,T]\times \mathbb{R}^{d}$\newline
\hspace*{0.2cm} (i) $\varphi (e^{\mu t}h^{+}(t,x))\in \mathbf{L}^{2}(%
\mathbb{R}^{d};\rho (x)dx)$,\newline \hspace*{0.2cm} (ii) $\left|
h(t,x)\right| \leq \kappa (1+\left| x\right| ^{\beta })$,\newline
here $h^{+}$ is the positive part of $h$.
\end{assumption}
\begin{assumption}\label{mark} $b:[0,T]\times \mathbb{R}^{d}\rightarrow %
\mathbb{R}^{d}$ and $\sigma :[0,T]\times \mathbb{R}^{d}\rightarrow \mathbb{R}%
^{d\times d}$ satisfy
\[
b\in C_{b}^{2}(\mathbb{R}^{d};\mathbb{R}^{d})\quad \mbox{and}\quad
\sigma \in C_{b}^{3}(\mathbb{R}^{d};\mathbb{R}^{d\times d}).
\]
\end{assumption}
We first study the following PDE

\[
\left\{
\begin{split}
& (\partial _{t}+\mathcal{L})\,u\;+\;F(t,x,u,\nabla u)=0,\quad \forall
\,(t,x)\,\in \,[0,T]\times \mathbb{R}^{d} \\
& u(x,T)=g(x),\quad \forall \,x\,\in \,\mathbb{R}^{d}
\end{split}
\right.
\]
where $F:[0,T]\times \mathbb{R}^{d}\times \mathbb{R}^{n}\times \mathbb{R}%
^{n\times d}\rightarrow \mathbb{R}$, such that
\[
F(t,x,u,p)=f(t,x,u,\sigma ^{*}p)
\]
and
\[
\mathcal{L}=\sum_{i=1}^{d}b_{i}\frac{\partial }{\partial x_{i}}+\frac{1}{2}%
\sum_{i,j=1}^{d}a_{i,j}\frac{\partial ^{2}}{\partial x_{i}\partial x_{j}},
\]
$a:=\sigma \sigma ^{*}$. Here $\sigma ^{*}$ is the transposed matrix of $%
\sigma $.

In order to study the weak solution of the PDE, we introduce the following
space
\[
\mathcal{H} := \{ u\in \mathbf{L}^{2}([0,T]\times \mathbb{R}^{d},ds\otimes
\rho (x)dx) \; \big|\; \sigma ^{*}\nabla u\in \mathbf{L} ^{2}(([0,T]\times %
\mathbb{R}^{d},ds\otimes \rho (x)dx)\}
\]
endowed with the norm
\[
\left\| u\right\| ^{2}:=\int_{\mathbb{R}^{d}}\int_{0}^{T}[\left|
u(s,x)\right| ^{2}+\left| (\sigma ^{*}\nabla u)(s,x)\right| ^{2}]\rho
(x)dsdx.
\]

\begin{definition}
\label{pde}We say that $u\in \mathcal{H}$ is the weak solution of the PDE
associated to $(g,f)$, if

(i) $\left\| u\right\| ^{2}<\infty ,$

(ii) for every $\phi \in C_{c}^{1,\infty }([0,T]\times \mathbb{R}^{d})$%
\begin{equation}
\int_{t}^{T}(u_{s},\partial _{t}\phi )ds+(u(t,\cdot ),\phi (t,\cdot
))-(g(\cdot ),\phi (\cdot ,T))+\int_{t}^{T}\mathcal{E}(u_{s},\phi
_{s})ds=\int_{t}^{T}(f(s,\cdot ,u_{s},\sigma ^{*}\nabla u_{s}),\phi _{s})ds.
\label{Pde1}
\end{equation}
where $(\phi ,\psi )=\int_{\mathbb{R}^{d}}\phi (x)\psi (x)dx$ denotes the
scalar product in $\mathbf{L}^{2}(\mathbb{R}^{d},dx)$ and
\[
\mathcal{E}(\psi ,\phi )=\int_{\mathbb{R}^{d}}((\sigma ^{*}\nabla \psi
)(\sigma ^{*}\nabla \phi )+\phi \nabla ((\frac{1}{2}\sigma ^{*}\nabla \sigma
+b)\psi ))dx
\]
is the energy of the system of our PDE which corresponds to the Dirichlet
form associated to the operator $\mathcal{L}$ when it is symmetric. Indeed $%
\mathcal{E}(\psi ,\phi )=-(\phi ,\mathcal{L}\psi )$.
\end{definition}

The probabilistic interpretation of the solution of PDE associated
with $g,f$, which satisfy Assumption
\ref{term-pde}-\ref{coef-pde1} was firstly studied by (Pardoux
\cite{Pardoux99}), where the author proved the existence of a
viscosity solution to this PDE, and gave its probabilistic
interpretation. In section 4, we consider the weak solution to PDE
(\ref{Pde1}) in Sobolev space, and give the proof of the existence
and uniqueness of the solution as well as the probabilistic
interpretation.

In the second part of this article, we will consider the obstacle
problem associated to the PDE (\ref{Pde1}) with obstacle function
$h$, where we restrict our study in the one dimensional case
($n=1$). Formulaly, The solution $u$ is
dominated by $h$, and verifies the equation in the following sense : $%
\forall (t,x)\,\in \,[0,T]\times \mathbb{R}^{d}$
\[
\begin{split}
& \text{(i) }(\partial _{t}+\mathcal{L})u+F(t,x,u,\nabla u)\leq 0,\quad \;%
\text{ on }\;\quad u(t,x)\geq h(t,x), \\
& \text{(ii) }(\partial _{t}+\mathcal{L})u+F(t,x,u,\nabla u)=0,\;\quad \text{
on }\quad \;u(t,x)> h(t,x), \\
& \text{(iii) }u(x,T)=g(x)\,.
\end{split}
\]
where $\mathcal{L}=\sum_{i=1}^{d}b_{i}\frac{\partial }{\partial x_{i}}+\frac{%
1}{2}\sum_{i,j=1}^{d}a_{i,j}\frac{\partial ^{2}}{\partial x_{i}\partial x_{j}%
}$, $a=\sigma \sigma ^{*}$. In fact, we give the following
formulation of the PDE with obstacle.

\begin{definition}
\label{o-pde}We say that $(u,\nu )$ is the weak solution of the PDE with
obstacle associated to $(g,f,h)$, if

(i) $\left\| u\right\| ^{2}<\infty $, $u\geq h$, and $u(T,x)=g(x)$.

(ii) $\nu $ is a positive Radon measure such that $\int_{0}^{T}\int_{\mathbb{R}%
^{d}}\rho (x)d\nu (t,x)<\infty ,$

(iii) for every $\phi \in C_{c}^{1,\infty }([0,T]\times \mathbb{R}^{d})$%
\begin{eqnarray}
&&\int_{t}^{T}(u_{s},\partial _{s}\phi )ds+(u(t,\cdot ),\phi (t,\cdot
))-(g(\cdot ),\phi (\cdot ,T))+\int_{t}^{T}\mathcal{E}(u_{s},\phi _{s})ds
\label{OPDE} \\
&=&\int_{t}^{T}(f(s,\cdot ,u_{s},\sigma ^{*}\nabla u_{s}),\phi
_{s})ds+\int_{t}^{T}\int_{\mathbb{R}^{d}}\phi (s,x)1_{\{u=h\}}d\nu (x,s).
\nonumber
\end{eqnarray}
\end{definition}

\section{Stochastic flow and random test functions}
\sectionmark{Stochastic flow and random test functions}

Let $(X_{s}^{t,x})_{t\leq s\leq T}$ be the solution of
\[
\left\{
\begin{split}
& dX_{s}^{t,x}=b(s,X_{s}^{t,x})ds+\sigma (s,X_{s}^{t,x})dB_{s}, \\
& X_{t}^{t,x}=x,
\end{split}
\right.
\]
where $b:[0,T]\times \mathbb{R}^{d}\rightarrow \mathbb{R}^{d}$ and
$\sigma :[0,T]\times \mathbb{R}^{d}\rightarrow \mathbb{R}^{d\times
d}$ satisfy Assumption \ref{mark}.

So $\{X_{s}^{t,x},x\in \mathbb{R}^d,t\leq s\leq T\}$ is the
stochastic flow associated to the diffuse $\{X_{s}^{t,x}\}$ and
denote by $\{\widehat{X}_{s}^{t,x},t\leq s\leq T\}$ the inverse
flow. It is known that $x\rightarrow \widehat{X}_{s}^{t,x}$ is
differentiable (Ikeda and Watanabe \cite{IW}). We denote by
$J(X_{s}^{t,x})$ the determinant of the Jacobian matrix of
$\widehat{X}_{s}^{t,x}$, which is positive, and
$J(X_{t}^{t,x})=1$.

For $\phi \in C_{c}^{\infty }(\mathbb{R}^{d})$ we define a process $\phi
_{t}:\Omega \times [0,T]\times \mathbb{R}^{d}\rightarrow \mathbb{R}$ by
\[
\phi _{t}(s,x):=\phi (\widehat{X}_{s}^{t,x})J(\widehat{X}_{s}^{t,x}).
\]
Following Kunita (See \cite{K82}), we know that for $v\in
\mathbf{L}^{2}(\mathbb{R}^{d})$, the composition of $v$ with the
stochastic flow is
\[
(v\circ X_{s}^{t,\cdot },\phi ):=(v,\phi _{t}(s,\cdot )).
\]
Indeed, by a change of variable, we have
\[
(v\circ X_{s}^{t,\cdot },\phi )=\int_{\mathbb{R}^{d}}v(y)\phi (\widehat{X}%
_{s}^{t,y})J(\widehat{X}_{s}^{t,y})dy=\int_{\mathbb{R}^{d}}v(X_{s}^{t,x})%
\phi (x)dx.
\]

The main idea in Bally and Matoussi \cite{BM2001} and Bally et al.
\cite{BCEF}, is to use $\phi _{t}$ as
a test function in (\ref{Pde1}) and (\ref{OPDE}). The problem is that $%
s\rightarrow \phi _{t}(s,x)$ is not differentiable so that $%
\int_{t}^{T}(u_{s},\partial _{s}\phi )ds$ has no sense. However
$\phi _{t}(s,x)$ is a semimartingale and they proved the following
semimartingale decomposition of $\phi _{t}(s,x)$:

\begin{lemma}
\label{comp}For every function $\phi \in \mathbf{C}_{c}^{2}(\mathbb{R}^{d})$%
,
\begin{equation}
\phi _{t}(s,x)=\phi (x)-\sum_{j=1}^{d}\int_{t}^{s}\left( \sum_{i=1}^{d}\frac{%
\partial }{\partial x_{i}}\left( \sigma _{ij}(x)\phi _{t}(r,x)\right)
\right) dB_{r}^{j}+\int_{t}^{s}\mathcal{L}^{*}\phi _{t}(r,x)dr,  \label{phi}
\end{equation}
where $\mathcal{L}^{*}$ is the adjoint operator of $\mathcal{L}$. So
\begin{equation}
d\phi _{t}(r,x)=-\sum_{j=1}^{d}\left( \sum_{i=1}^{d}\frac{\partial }{%
\partial x_{i}}\left( \sigma _{ij}(x)\phi _{t}(r,x)\right) \right)
dB_{r}^{j}+\mathcal{L}^{*}\phi _{t}(r,x)dr,  \label{dphi}
\end{equation}
\end{lemma}

Then in (\ref{Pde1}) we may replace $\partial _{s}\phi ds$ by the It\^{o}
stochastic integral with respect to $d\phi _{t}(s,x)$, and have the
following proposition which allows us to use $\phi _{t}$ as a test function.
The proof will be given in the appendix.

\begin{proposition}
\label{test1} Assume that assumptions \ref{term-pde},
\ref{coef-pde} and \ref {coef-pde1} hold. Let $u\in \mathcal{H}$
be a weak solution of PDE (\ref {Pde1}), then for $s\in [t,T]$ and
$\phi \in C_{c}^{2}(\mathbb{R}^{d})$,
\begin{eqnarray}
&&\ \int_{\mathbb{R}^{d}}\int_{s}^{T}u(r,x)d\phi _{t}(r,x)dx-(g(\cdot ),\phi
_{t}(T,\cdot ))+(u(s,\cdot ),\phi _{t}(s,\cdot ))-\int_{s}^{T}\mathcal{E}%
(u(r,\cdot ),\phi _{t}(r,\cdot ))dr  \nonumber \\
\  &=&\int_{\mathbb{R}^{d}}\int_{s}^{T}f(r,x,u(r,x),\sigma ^{*}\nabla
u(r,x))\phi _{t}(r,x)drdx.\text{ a.s.}  \label{tesf2}
\end{eqnarray}
\end{proposition}

\begin{remark}
Here $\phi _{t}(r,x)$ is $\mathbb{R}$-valued. We consider that in
(\ref{tesf2}), the equality holds for each component of $u$.
\end{remark}

We need the result of equivalence of norms, which play important
roles in existence proof for PDE under monotonic conditions. The
equivalence of functional norm and stochastic norm is first proved
by Barles and Lesigne \cite{BL} for $\rho =1$. In Bally and
Matoussi \cite{BM2001} proved the same result for weighted
integrable function by using probabilistic method. Let $\rho $ be
a weighted function, we take $\rho (x):=\exp (F(x))$, where
$F:\mathbb{R}^{d}\rightarrow \mathbb{R}$
is a continuous function. Moreover, we assume that there exists a constant $%
R>0$, such that for $\left| x\right| >R$, $F\in C_{b}^{2}(\mathbb{R}^{d},%
\mathbb{R})$. For instant, we can take $\rho (x)=(1+\left|
x\right| )^{-q}$
or $\rho (x)=\exp \alpha \left| x\right| $, with $q>d+1$, $\alpha \in \mathbb{R}%
. $

\begin{proposition}
\label{equiv}Suppose that assumption \ref{mark} hold, then there
exists
two constants $k_{1}$, $k_{2}>0$, such that for every $t\leq s\leq T$ and $%
\phi \in \mathbf{L}^{1}(\mathbb{R}^{d},\rho (x)dx)$, we have
\begin{equation}
k_{2}\int_{\mathbb{R}^{d}}\left| \phi (x)\right| \rho (x)dx\leq \int_{%
\mathbb{R}^{d}}E(\left| \phi (X_{s}^{t,x})\right| )\rho (x)dx\leq k_{1}\int_{%
\mathbb{R}^{d}}\left| \phi (x)\right| \rho (x)dx,  \label{equiv1}
\end{equation}
Moreover, for every $\psi \in \mathbf{L}^{1}([0,T]\times \mathbb{R}%
^{d},dt\otimes \rho (x)dx)$%
\begin{eqnarray}
k_{2}\int_{\mathbb{R}^{d}}\int_{t}^{T}\left| \psi (s,x)\right| \rho (x)dsdx
&\leq &\int_{\mathbb{R}^{d}}\int_{t}^{T}E(\left| \psi (s,X_{s}^{t,x})\right|
)\rho (x)dsdx  \label{equiv2} \\
\  &\leq &k_{1}\int_{\mathbb{R}^{d}}\int_{t}^{T}\left| \psi (s,x)\right|
\rho (x)dsdx,  \nonumber
\end{eqnarray}
where the constants $k_{1}$, $k_{2}$ depend only on $T$, $\rho $ and the
bounds of the first (resp. first and second) derivatives of $b$ (resp. $%
\sigma $).
\end{proposition}

This proposition is easy to get from the follwing Lemma, see Lemma
5.1 in Bally and Matoussi \cite{BM2001}.

\begin{lemma}
\label{flo-bon}There exist two constants $c_{1}>0$ and $c_{2}>0$ such that $%
\forall x\in \mathbb{R}^{d}$, $0\leq t\leq T$%
\[
c_{1}\leq E\left( \frac{\rho (t,\widehat{X}_{t}^{0,x})J(\widehat{X}%
_{t}^{0,x})}{\rho (x)}\right) \leq c_{2}.
\]
\end{lemma}

\section{Sobolev's Solutions for PDE's under monotonicity condition}
\sectionmark{Sobolev Solutions for PDE}

In this section we shall  study the solution of the PDE whose
coefficient $f$ satisfies the monotonicity condition. For this sake,
we introduce the BSDE associated with $(g,f)$: for $t\leq s\leq T$,
\begin{equation}
Y_{s}^{t,x}=g(X_{T}^{t,x})+%
\int_{s}^{T}f(r,X_{r}^{t,x},Y_{r}^{t,x},Z_{r}^{t,x})dr-%
\int_{s}^{T}Z_{s}^{t,x}dB_{s}.  \label{bsde1}
\end{equation}
Thanks to  the equivalence of the norms result (\ref{equiv}), we
know that  $g(X_{T}^{t,x})$ and $f(s,X_{s}^{t,x},0,0)$ make sense in
the BSDE (\ref{bsde1}). Moreover we have \[ g(X_{T}^{t,x}) \in
\mathbf{L}_n^2({\cal F}_T) \mbox{ and } f(.,X_{.}^{t,x},0,0) \in
\mathbf{H}_n^2(0,T).
\]
It follows from the results from Pardoux \cite{Pardoux99} that for each $%
(t,x)$, there exists a unique pair $(Y^{t,x},Z^{t,x})\in
\mathbf{S}^{2}(t,T)\times \mathbf{H}_{n\times d}^{2}(t,T)$ of
$\{\mathcal{F}_{s}^{t}\}$ progressively measurable processes, which
solves this BSDE$(g,f)$. \\[0.1cm] The main result of this section is

\begin{theorem}
\label{mr1}Suppose that assumptions \ref{term-pde}-\ref{coef-pde1} and \ref
{bar-pde} hold. Then there exists a unique weak solution $u\in \mathcal{%
H}$ of the PDE (\ref{Pde1}).  Moreover we have the probabilistic
interpretation of the solution:
\begin{equation}
\label{F-K} u(t,x)=Y_{t}^{t,x},\quad  (\sigma^{*}\nabla
u)(t,x)=Z_{t}^{t,x}, \quad dt\otimes dx -a.e.
\end{equation}
and moreover $ \; \; Y_{s}^{t,x}=u(s,X_{s}^{t,x})$,
$Z_{s}^{t,x}=(\sigma ^{*}\nabla u)(s,X_{s}^{t,x})$, $dt \otimes dP
\otimes dx $-a.e. $\forall s \in [t, T]$.
\end{theorem}
\textbf{Proof}: We start to prove the existence result. \\[0.2cm]
\textbf{\large a) Existence} :  We prove the existence in three
steps. By integration by parts formula, we know that $u$ solves
(\ref{Pde1}) if and only if
\begin{equation*}
\widehat{u}(t,x)=e^{\mu t}u(t,x)  \label{tran-u}
\end{equation*}
is a solution of the PDE$(\widehat{g},\widehat{f})$, where
\begin{equation}
\label{tran-m}
 \widehat{g}(x) =e^{\mu T}g(x) \;\; \text{and}
  \; \; \widehat{f}(t,x,y,z) =e^{\mu t}f(t,x,e^{-\mu t}y,e^{-\mu t}z)-\mu y.
\end{equation}
 Then the coefficient $\widehat{f}$ satisfies the
assumption \ref{coef-pde1} as $f$, except that \ref{coef-pde1}-(iii)
is replaced by

\begin{equation}
(y-y^{\prime })(f(t,x,y,z)-f(t,x,y^{\prime },z))\leq 0.  \label{tran-f}
\end{equation}

In the first two steps, we consider the case where $f$ does not depend on $%
\nabla u$, and write $f(t,x,y)$ for $f(t,x,y,v(t,x))$, where $v$ is in $%
\mathbf{L}^{2}([0,T]\times \mathbb{R}^{d},dt\otimes \rho
(x)dx)$.\\[0.1cm]

We assume first that $f(t,x,y)$ satisfies the following  assumption
 \ref{coef-pde1}\textbf{':} $%
\forall (t,x,y,y^{\prime })\in [0,T]\times \mathbb{R}^{d}\times \mathbb{R}%
^{n}\times \mathbb{R}^{n}$,

\begin{enumerate}
\item[(i)]  $\left| f(t,x,y)\right| \leq \left| f(t,x,0)\right| +\varphi
(\left| y\right| )$,

\item[(ii)]  $\left\langle y-y^{\prime },f(t,x,y)-f(t,x,y^{\prime
})\right\rangle \leq 0$,

\item[(iii)]  $y\rightarrow f(t,x,y)$ is continuous, $\forall (t,x)\in
[0,T]\times \mathbb{R}^{d}$.
\end{enumerate}

\textbf{\large Step 1} : \textit{ Suppose that $g(x)$, $f(t,x,0)$
are uniformly bounded}, i.e. there exists a constant $C$, such that
\begin{equation}
\left| g(x)\right| +\sup_{0\leq t\leq T}\left| f(t,x,0)\right| \leq
C \label{con-b}
\end{equation}
where $C$ as a constant which can be changed line by line.

Define $ f_{n}(t,y):=(\theta _{n}*f(t,\cdot ))(y)$  where $\theta
_{n}:\mathbb{R}^{n}\rightarrow \mathbb{R}_{+}$ is a sequence of
smooth functions with compact support, which approximate the Dirac
distribution at $0$, and satisfy $\int \theta _{n}(z)dz=1$. Let
 $ \{(Y_{s}^{n,t,x},Z_{s}^{n,t,x}), \; t \leq s\leq T \}$ be  the solution
 of BSDE associated to $(g(X_{T}^{t,x}),f_{n})$, namely,
\begin{equation}
Y_{s}^{n,t,x}=g(X_{T}^{t,x})+\int_{s}^{T}f(r,X_{r}^{t,x},Y_{r}^{n,t,x})dr-%
\int_{s}^{T}Z_{r}^{n,t,x}dB_{r},\; \text{P-a.s..}  \label{bsde-bn}
\end{equation}
Then for each $n \in \mathbb{N}$, we have
\[
\left| Y_{s}^{n,t,x}\right| \leq e^{T}C,
\]
and
\[
\left| f_{n}(s,X_{s}^{t,x},Y_{s}^{n,t,x})\right| ^{2}\leq 2\left|
f_{n}(s,X_{s}^{t,x},0)\right| ^{2}+2\psi ^{2}(e^{\frac{T}{2}}\sqrt{C})
\]
where $\psi (r):=\sup_{n}\sup_{\left| y\right| \leq r}\int_{\mathbb{R}%
^{n}}\varphi (\left| y\right| )\theta _{n}(y-z)dz$. So there exists a
constant $C>0$, s.t.
\begin{equation}
\sup_{n}\int_{\mathbb{R}^{d}}E\int_{t}^{T}(\left| Y_{s}^{n,t,x}\right|
^{2}+\left| f_{n}(s,X_{s}^{t,x},Y_{s}^{n,t,x})\right| ^{2}+\left|
Z_{s}^{n,t,x}\right| ^{2})\rho (x)dsdx\leq C.  \label{est-n}
\end{equation}
Then let $n\rightarrow \infty $ on the both sides of (\ref{bsde-bn}), we get
that the  limit $(Y_{s}^{t,x},Z_{s}^{t,x})$ of $%
(Y_{s}^{n,t,x},Z_{s}^{n,t,x})$, satisfies
\begin{equation}
Y_{s}^{t,x}=g(X_{T}^{t,x})+\int_{s}^{T}f(r,X_{r}^{t,x},Y_{r}^{t,x})dr-%
\int_{s}^{T}Z_{r}^{t,x}dB_{r},\; \text{P-a.s..}  \label{bsde-b}
\end{equation}
Moreover we obtain from the estimate (\ref{est-n}) that
\begin{equation}
\label{est-BSDE}
 \int_{\mathbb{R}^{d}}\int_{t}^{T}E(\left|
Y_{s}^{t,x}\right| ^{2}+\left| Z_{s}^{t,x}\right| ^{2})\rho
(x)dsdx<\infty .
\end{equation}
Notice that $(Y_{t}^{t,x},Z_{t}^{t,x})$
are ${\cal F}^t_t$ measurable, which implies they are deterministic.
Define $u(t,x):=Y_{t}^{t,x}$, and $v(t,x):=Z_{t}^{t,x}$.
 By the flow
property of  $X_{r}^{s,x}$ and by the uniqueness of the solution of the BSDE (\ref{bsde-b}),
we have that $%
Y_{s}^{t,x}=u(s,X_{s}^{t,x})$ and $Z_{s}^{t,x}=v(s,X_{s}^{t,x})$.\\
The terminal condition $g$ and $ f(.,.,0,0)$ are not continuous in
$t$ and $x$, and assumed to  belong in  a suitable weighted $L^2$
space, so the solution $u$ and for instance $v$ are not in general
continuous, and are only defined a.e. in $ [0,T] \times
\mathbb{R}^d$. So in order to give meaning to the expression $u (s,
X_s^{t,x})$ (resp. $v (s, X_s^{t,x})$), and following Bally and
Matoussi \cite{BM2001}, we apply a regularization procedure on the
final condition $g$ and the coefficient $f$. Actually, according to
Pardoux and Peng (\cite{PP92}, Theorem 3.2), if the  coefficient $
(g, f )$ are smooth, then the PDE (\ref{Pde1}) admits a unique
classical solution $u \in C^{1,2} ([0,T]\times \mathbb{R}^d)$.
Therefore  the approximated expression $u (s, X_s^{t,x})$ (resp. $v
(s, X_s^{t,x})$) has a meaning   and then pass to the limit in $L^2$
spaces like us in Bally and Matoussi \cite{BM2001}. \\
Now, the equivalence of norm result (\ref{equiv2}) and estimate
(\ref{est-BSDE})  follow that
 $u,v\in \mathbf{L}^{2}([0,T]\times \mathbb{R}^{d},dt\otimes \rho (x)dx)$.
 Finally, let $F(r,x)=f(r,X_{r}^{t,x},Y_{r}^{t,x})$, we know that $F(s,x)\in \mathbf{%
L}^{2}([0,T]\times \mathbb{R}^{d},dt\otimes \rho (x)dx)$, in view of
\begin{eqnarray*}
\int_{\mathbb{R}^{d}}\int_{t}^{T}\left| F(s,x)\right| ^{2}\rho (x)dsdx &\leq
&\frac{1}{k_{2}}\int_{\mathbb{R}^{d}}\int_{t}^{T}E\left|
F(s,X_{s}^{t,x})\right| ^{2}\rho (x)dsdx \\
&=&\frac{1}{k_{2}}\int_{\mathbb{R}^{d}}\int_{t}^{T}E\left|
f(s,X_{s}^{t,x},Y_{s}^{t,x})\right| ^{2}\rho (x)dsdx<\infty .
\end{eqnarray*}
So that from theorem 2.1 in \cite{BM2001}, we get that $v=\sigma
^{*}\nabla u$ and that $u\in \mathcal{H}$ solves the PDE associated
to $(g,f)$ under the bounded assumption.\\[0.2cm]
\textbf{\large Step 2} :  \textit{We assume $g\in
\mathbf{L}^{2}(\mathbb{R}^{d},\rho (x)dx)$, $f$ satisfies the
assumption \ref{coef-pde1}\textbf{'} and $f(t,x,0)\in
\mathbf{L}^{2}([0,T]\times \mathbb{R}^{d},dt\otimes \rho (x)dx)$}.
We approximate $g$ and $f$ by bounded functions as follows :
\begin{eqnarray}
g_{n}(x) &=&\Pi _{n}(g(x)),  \label{appr} \\
f_{n}(t,x,y) &=&f(t,x,y)-f(t,x,0)+\Pi _{n}(f(t,x,0)),  \nonumber
\end{eqnarray}
where
\[
\Pi _{n}(y):=\frac{\min (n,\left| y\right| )}{\left| y\right| }y.
\]
Clearly, the pair $(g_{n},f_{n})$ satisfies the assumption (\ref{con-b}) of
step 1, and
\begin{eqnarray}
g_{n} &\rightarrow &g\text{ in }\mathbf{L}^{2}(\mathbb{R}^{d},\rho (x)dx),
\label{convre} \\
f_{n}(t,x,0) &\rightarrow &f(t,x,0)\text{ in }\mathbf{L}^{2}([0,T]\times %
\mathbb{R}^{d},dt\otimes \rho (x)dx).  \nonumber
\end{eqnarray}

Denote $(Y_{s}^{n,t,x},Z_{s}^{n,t,x})\in
\mathbf{S}_{n}^{2}(t,T)\times \mathbf{H}_{n\times d}^{2}(t,T)$ the
solution of the BSDE$(\xi _{n},f_{n})$, where $\xi
_{n}=g_{n}(X_{T}^{t,x})$, i.e.
\[
Y_{s}^{n,t,x}=g_{n}(X_{T}^{t,x})+%
\int_{s}^{T}f_{n}(r,X_{r}^{t,x},Y_{r}^{n,t,x})dr-%
\int_{s}^{T}Z_{r}^{n,t,x}dB_{r}.
\]
Then from the results in step 1, $u_{n}(t,x)=Y_{t}^{n,t,x}$ and $%
u_{n}(t,x)\in \mathcal{H}$, is the weak solution of the PDE$(g_{n},f_{n})$,
with
\begin{equation}
Y_{s}^{n,t,x}=u_{n}(s,X_{s}^{t,x}),Z_{s}^{n,t,x}=(\sigma ^{*}\nabla
u_{n})(s,X_{s}^{t,x}),\text{a.s.}  \label{relation4}
\end{equation}

For $m,n\in \mathbb{N}$, applying It\^{o}'s formula to $\left|
Y_{s}^{m,t,x}-Y_{s}^{n,t,x}\right| ^{2}$, we get
\begin{equation}
\label{Ito1}
\begin{split}
 E\left| Y_{s}^{m,t,x}-Y_{s}^{n,t,x}\right|
^{2}&+E\int_{s}^{T}\left| Z_{r}^{m,t,x}-Z_{r}^{n,t,x}\right| ^{2}dr
\leq E\left| g_{m}(X_{T}^{t,x})-g_{n}(X_{T}^{t,x})\right| ^{2}\\
&+E\int_{s}^{T}\left| Y_{r}^{m,t,x}-Y_{r}^{n,t,x}\right|
^{2}dr+E\int_{s}^{T}\left|
f_{m}(r,X_{r}^{t,x},0)-f_{n}(r,X_{r}^{t,x},0)\right| ^{2}dr.
\end{split}
\end{equation} From the equivalence of the norms (\ref{equiv1}) and
(\ref{equiv2}), it follows
\begin{equation*}
\begin{split}
\int_{\mathbb{R}^{d}} & E\left| Y_{s}^{m,t,x}-Y_{s}^{n,t,x}\right|
^{2}\rho (x)dx \leq \int_{\mathbb{R}^{d}}E\left|
g_{m}(X_{T}^{t,x})-g_{n}(X_{T}^{t,x})\right| ^{2}\rho (x)dx\\
 & +\int_{\mathbb{R}%
^{d}}E\int_{s}^{T}\left| Y_{r}^{m,t,x}-Y_{r}^{n,t,x}\right|
^{2}dr\rho (x)dx  +\int_{\mathbb{R}^{d}}E\int_{s}^{T}\left|
f_{m}(r,X_{r}^{t,x},0)-f_{n}(r,X_{r}^{t,x},0)\right| ^{2}dr\rho (x)dx \\
&\leq \int_{\mathbb{R}^{d}}E\int_{s}^{T}\left|
Y_{r}^{m,t,x}-Y_{r}^{n,t,x}\right| ^{2}dr\rho (x)dx+k_{1}\int_{\mathbb{R}%
^{d}}E\left| g_{m}(x)-g_{n}(x)\right| ^{2}\rho (x)dx \\
&\quad \quad \quad +k_{1}\int_{\mathbb{R}^{d}}\int_{t}^{T}\left|
f_{m}(r,x,0)-f_{n}(r,x,0)\right| ^{2}\rho (x)drdx,
\end{split}
\end{equation*}
and by Gronwall's inequality and (\ref{convre}), we get as $m,n\rightarrow
\infty $%
\[
\sup_{t\leq s\leq T}\int_{\mathbb{R}^{d}}E\left|
Y_{s}^{m,t,x}-Y_{s}^{n,t,x}\right| ^{2}\rho (x)dx\rightarrow 0.
\]
It follows immediately as $m,n\rightarrow \infty $%
\[
\int_{\mathbb{R}^{d}}E\int_{s}^{T}\left| Y_{r}^{m,t,x}-Y_{r}^{n,t,x}\right|
^{2}\rho (x)drdx+\int_{\mathbb{R}^{d}}E\int_{s}^{T}\left|
Z_{r}^{m,t,x}-Z_{r}^{n,t,x}\right| ^{2}\rho (x)drdx\rightarrow 0.
\]
Using again the equivalence of the norms (\ref{equiv2}), we get:
\begin{eqnarray*}
&&\int_{t}^{T}\int_{\mathbb{R}^{d}}\left| u_{m}(s,x)-u_{n}(s,x)\right|
^{2}+\left| \sigma ^{*}\nabla u_{m}(s,x)-\sigma ^{*}\nabla u_{n}(s,x)\right|
^{2}\rho (x)dxds \\
&\leq &\frac{1}{k_{2}}\int_{t}^{T}\int_{\mathbb{R}^{d}}E(\left|
u_{m}(s,X_{s}^{t,x})-u_{n}(s,X_{s}^{t,x})\right| ^{2}+\left| \sigma
^{*}\nabla u_{m}(s,X_{s}^{t,x})-\sigma ^{*}\nabla
u_{n}(s,X_{s}^{t,x})\right| ^{2})\rho (x)dsdx \\
&=&\frac{1}{k_{2}}\int_{t}^{T}\int_{\mathbb{R}^{d}}E(\left|
Y_{s}^{m,t,x}-Y_{s}^{n,t,x}\right| ^{2}+\left|
Z_{s}^{m,t,x}-Z_{s}^{n,t,x}\right| ^{2})\rho (x)dsdx\rightarrow 0.
\end{eqnarray*}
as $m,n\rightarrow \infty $, i.e. $\{u_{n}\}$ is Cauchy sequence in $%
\mathcal{H}$. Denote its limit as $u$, so $u\in \mathcal{H}$, and satisfies
for every $\phi \in C_{c}^{1,\infty }([0,T]\times \mathbb{R}^{d})$,
\begin{equation}
\int_{t}^{T}(u_{s},\partial _{t}\phi )ds+(u(t,\cdot ),\phi (t,\cdot
))-(g(\cdot ),\phi (\cdot ,T))+\int_{t}^{T}\mathcal{E}(u_{s},\phi
_{s})ds=\int_{t}^{T}(f(s,\cdot ,u_{s}),\phi _{s})ds.  \label{Pde3}
\end{equation}
On the other hand, $(Y_{\cdot }^{n,t,x},Z_{\cdot }^{n,t,x})$ converges to $%
(Y_{\cdot }^{t,x},Z_{\cdot }^{t,x})$ in $\mathbf{S}_{n}^{2}(0,T)\times
\mathbf{H}_{n\times d}^{2}(0,T)$, which is the solution of the BSDE with
parameters $(g(X_{T}^{t,x}),f)$; by the equivalence of the norms, we deduce
that
\[
Y_{s}^{t,x}=u(s,X_{s}^{t,x}),Z_{s}^{t,x}=\sigma ^{*}\nabla u(s,X_{s}^{t,x}),%
\text{ a.s. }\forall s\in [t,T],
\]
specially $Y_{t}^{t,x}=u(t,x)$, $Z_{t}^{t,x}=\sigma ^{*}\nabla
u(t,x)$.

Now, it's easy to the generalize the result to the case when $f$
satisfies
assumption \ref{coef-pde} .\\[0.2cm]

\textbf{\large Step 3}: \textit{In this step, we consider the case where $f$ depends on $%
\nabla u$}. Assume that $g$, $f$ satisfy the assumptions \ref{term-pde} - \ref{coef-pde1}%
, with assumption \ref{coef-pde1}-(iii) replaced by (\ref{tran-f}). From the
result in step 2, for any given $n\times d$-matrix-valued function $v\in
\mathbf{L}^{2}([0,T]\times \mathbb{R}^{d},dt\otimes \rho (x)dx)$, $%
f(t,x,u,v(t,x))$ satisfies the assumptions in step 2. So the PDE$%
(g,f(t,x,u,v(t,x)))$ admits a unique solution $u\in \mathcal{H}$ satisfying
(i) and (ii) in the definition \ref{pde}.

Set $V_{s}^{t,x}=v(s,X_{s}^{t,x})$, then $V_{s}^{t,x}\in \mathbf{H}_{n\times
d}^{2}(0,T)$ in view of the equivalence of the norms. We consider the
following BSDE with solution $(Y_{\cdot }^{t,x},Z_{\cdot }^{t,x})$%
\[
Y_{s}^{t,x}=g(X_{T}^{t,x})+%
\int_{s}^{T}f(s,X_{s}^{t,x},Y_{s}^{t,x},V_{s}^{t,x})ds-%
\int_{s}^{T}Z_{s}^{t,x}dB_{s},
\]
then $Y_{s}^{t,x}=u(s,X_{s}^{t,x})$, $Z_{s}^{t,x}=\sigma ^{*}\nabla
u(s,X_{s}^{t,x})$, a.s. $\forall s\in [t,T]$.

Now we can construct a mapping $\Psi $ from $\mathcal{H}$ into itself. For
any $\overline{u}\in \mathcal{H}$, $u=\Psi (\overline{u})$ is the weak
solution of the PDE with parameters $g(x)$ and $f(t,x,u,\sigma ^{*}\nabla
\overline{u})$.

Symmetrically we introduce a mapping $\Phi $ from $\mathbf{H}_{n%
}^{2}(t,T)\times \mathbf{H}_{n\times d}^{2}(t,T)$ into itself. For any $(U^{t,x},V^{t,x})\in \mathbf{H}_{%
n}^{2}(t,T)\times \mathbf{H}_{n\times d}^{2}(t,T)$,
$(Y^{t,x},Z^{t,x})=\Phi (U^{t,x},V^{t,x})$ is the
solution of the BSDE with parameters $g(X_{T}^{t,x})$ and $%
f(s,X_{s}^{t,x},Y_{s}^{t,x},V_{s}^{t,x})$. Set $V_{s}^{t,x}=\sigma
^{*}\nabla \overline{u}(s,X_{s}^{t,x})$, then $Y_{s}^{t,x}=u(s,X_{s}^{t,x})$%
, $Z_{s}^{t,x}=\sigma ^{*}\nabla u(s,X_{s}^{t,x})$, a.s.a.e..

Let $\overline{u}_{1}$, $\overline{u}_{2}\in \mathcal{H}$, and $u_{1}=\Psi (%
\overline{u}_{1})$, $u_{2}=\Psi (\overline{u}_{2})$, we consider
the
difference $\triangle u: =u_{1}-u_{2}$, $\triangle \overline{u}: =\overline{u}%
_{1}-\overline{u}_{2}$. Set $V_{s}^{t,x,1}:=\sigma ^{*}\nabla \overline{u}%
_{1}(s,X_{s}^{t,x})$, $V_{s}^{t,x,2}:=\sigma ^{*}\nabla \overline{u}%
_{2}(s,X_{s}^{t,x})$. We denote by $(Y^{t,x,1},Z^{t,x,1})$(resp. $%
(Y^{t,x,2},Z^{t,x,2})$) the solution of the BSDE with parameters $%
g(X_{T}^{t,x})$ and $f(s,X_{s}^{t,x},Y_{s}^{t,x},V_{s}^{t,x,1})$ (resp. $%
f(s,X_{s}^{t,x},Y_{s}^{t,x},V_{s}^{t,x,2})$); then for a.e. $\forall
s\in [t,T],$
\begin{eqnarray*}
Y_{s}^{t,x,1} &=&u_{1}(s,X_{s}^{t,x}),Z_{s}^{t,x,1}=\sigma
^{*}\nabla
u_{1}(s,X_{s}^{t,x}), \\
Y_{s}^{t,x,2} &=&u_{2}(s,X_{s}^{t,x}),Z_{s}^{t,x,2}=\sigma
^{*}\nabla u_{2}(s,X_{s}^{t,x}),
\end{eqnarray*}
Denote $\triangle Y_{s}^{t,x}:=Y_{s}^{t,x,1}-Y_{s}^{t,x,2}$,
$\triangle Z_{s}^{t,x}:=Z_{s}^{t,x,1}-Z_{s}^{t,x,2}$, $\triangle
V_{s}^{t,x}:=V_{s}^{t,x,1}-V_{s}^{t,x,2}$. By It\^{o}'s formula applied to $%
e^{\gamma t}E\left| \triangle Y_{s}^{t,x}\right| ^{2}$, for some $\alpha$ and $\gamma \in %
\mathbb{R}$, we have
\[
e^{\gamma t}E\left| \triangle Y_{s}^{t,x}\right|
^{2}+E\int_{s}^{T}e^{\gamma s}(\gamma \left| \triangle
Y_{r}^{t,x}\right| ^{2}+\left| \triangle Z_{r}^{t,x}\right|
^{2})dr\leq E\int_{s}^{T}e^{\gamma s}(\frac{k^{2}}{\alpha }\left|
\triangle Y_{r}^{t,x}\right| ^{2}+\alpha \left| \triangle
V_{r}^{t,x}\right| ^{2})dr,
\]
Using the equivalence of the norms, we deduce that
\begin{eqnarray*}
&&\int_{\mathbb{R}^{d}}\int_{t}^{T}e^{\gamma s}(\gamma \left|
\triangle u(s,x)\right| ^{2}+\left| \sigma ^{*}\nabla (\triangle
u)(s,x)\right|
^{2})\rho (x)dsdx \\
&\leq &\frac{1}{k_{2}}\int_{\mathbb{R}^{d}}\int_{t}^{T}e^{\gamma
s}E(\gamma \left| \triangle Y_{r}^{t,x}\right| ^{2}+\left| \triangle
Z_{r}^{t,x}\right|
^{2})\rho (x)drdx \\
&\leq &\frac{1}{k_{2}}\int_{\mathbb{R}^{d}}\int_{s}^{T}e^{\gamma s}E(\frac{%
k^{2}}{\alpha }\left| \triangle Y_{r}^{t,x}\right| ^{2}+\alpha
\left|
\triangle V_{r}^{t,x}\right| ^{2})\rho (x)drdx \\
&\leq &\frac{k_{1}}{k_{2}}\int_{\mathbb{R}^{d}}\int_{s}^{T}e^{\gamma s}(%
\frac{k^{2}}{\alpha }\left| \triangle u(s,x)\right| ^{2}+\alpha
\left| \sigma ^{*}\nabla (\triangle \overline{u})(s,x)\right|
^{2})\rho (x)dsdx.
\end{eqnarray*}
Set $\alpha =\frac{k_{2}}{2k_{1}}$, $\gamma =1+\frac{2k_{1}^{2}}{k_{2}^{2}}%
k^{2}$, then we get
\begin{eqnarray*}
&&\int_{\mathbb{R}^{d}}\int_{t}^{T}e^{\gamma s}(\left| \triangle
u(s,x)\right| ^{2}+\left| \sigma ^{*}\nabla (\triangle
u)(s,x)\right|
^{2})\rho (x)dsdx \\
&\leq &\frac{1}{2}\int_{\mathbb{R}^{d}}\int_{t}^{T}e^{\gamma
s}\left| \sigma
^{*}\nabla (\triangle \overline{u})(s,x)\right| ^{2}\rho (x)dsdx, \\
&\leq &\frac{1}{2}\int_{\mathbb{R}^{d}}\int_{t}^{T}e^{\gamma
s}(\left| \triangle \overline{u}(s,x)\right| ^{2}+\left| \sigma
^{*}\nabla (\triangle \overline{u})(s,x)\right| ^{2})\rho (x)dsdx.
\end{eqnarray*}
Consequently, $\Psi $ is a strict contraction on $\mathcal{H}$
equipped with the norm
\[
\left\| u\right\| _{\gamma
}^{2}:=\int_{\mathbb{R}^{d}}\int_{t}^{T}e^{\gamma s}(\left|
u(s,x)\right| ^{2}+\left| \sigma ^{*}\nabla u(s,x)\right|
^{2})\rho (x)dsdx.
\]
So $\Psi $ has  fixed point $u \in \mathcal{H} $ which is  the
solution  of the PDE (\ref{Pde1}) associated to $(g,f)$. Moreover,
for $t\leq s\leq T$,
\[
Y_{s}^{t,x}=u(s,X_{s}^{t,x}),Z_{s}^{t,x}=\sigma ^{*}\nabla u(s,X_{s}^{t,x}),%
\text{ .a.e.}
\]
and specially $Y_{t}^{t,x}=u(t,x),Z_{t}^{t,x}=\sigma ^{*}\nabla u(t,x), \, \, \text{a.e}.$ \\[%
0.2cm] \textbf{\large b) Uniqueness } :  Let $u^{1}$ and $u^{2}\in
\mathcal{H}$ be two solutions of the PDE$(g,f)$.
From Proposition \ref{test1}, for $\phi \in C_{c}^{2}(\mathbb{R}^{d})$ and $%
i=1,2$%
\begin{eqnarray}
&&\int_{\mathbb{R}^{d}}\int_{s}^{T}u^{i}(r,x)d\phi
_{t}(r,x)dx+(u^{i}(s,\cdot ),\phi _{t}(s,\cdot ))-(g(\cdot ),\phi _{t}(\cdot
,T))-\int_{s}^{T}\mathcal{E(}u^{i}(r,\cdot ),\phi _{t}(r,\cdot ))dr
\nonumber \\
&=&\int_{s}^{T}\int_{\mathbb{R}^{d}}\phi _{t}(r,x)f(r,x,u^{i}(r,x),\sigma
^{*}\nabla u^{i}(r,x))drdx.  \label{pde-u1}
\end{eqnarray}
By (\ref{phi}), we get
\begin{eqnarray*}
\int_{\mathbb{R}^{d}}\int_{s}^{T}u^{i}d\phi _{t}(r,x)dx
&=&\int_{s}^{T}(\int_{\mathbb{R}^{d}}(\sigma ^{*}\nabla u^{i})(r,x)\phi
_{t}(r,x)dx)dB_{r} \\
&&+\int_{s}^{T}\int_{\mathbb{R}^{d}}\left( (\sigma ^{*}\nabla u^{i})(\sigma
^{*}\nabla \phi _{r})+\phi \nabla ((\frac{1}{2}\sigma ^{*}\nabla \sigma
+b)u_{r}^{i})\right) dxdr.
\end{eqnarray*}

We substitute this in (\ref{pde-u1}), and get
\begin{eqnarray*}
\int_{\mathbb{R}^{d}}u^{i}(s,x)\phi _{t}(s,x)dx &=&(g(\cdot ),\phi
_{t}(\cdot ,T))-\int_{s}^{T}\int_{\mathbb{R}^{d}}(\sigma ^{*}\nabla
u^{i})(r,x)\phi _{t}(r,x)dxdB_{r} \\
&&+\int_{s}^{T}\int_{\mathbb{R}^{d}}\phi _{t}(r,x)f(r,x,u^{i}(r,x),\sigma
^{*}\nabla u^{i}(r,x))drdx.
\end{eqnarray*}
Then by the change of variable $y=\widehat{X}_{r}^{t,x}$, we obtain
\begin{eqnarray*}
\int_{\mathbb{R}^{d}}u^{i}(s,X_{s}^{t,y})\phi (y)dy &=&\int_{\mathbb{R}%
^{d}}g(X_{T}^{t,y})\phi (y)dy+\int_{s}^{T}\int_{\mathbb{R}^{d}}\phi
(y)f(s,X_{s}^{t,y},u^{i}(s,X_{s}^{t,y}),\sigma ^{*}\nabla
u^{i}(s,X_{s}^{t,y}))dyds \\
&&-\int_{s}^{T}\int_{\mathbb{R}^{d}}(\sigma ^{*}\nabla
u^{i})(r,X_{r}^{t,y})\phi (y)dydB_{r}.
\end{eqnarray*}
Since $\phi $ is arbitrary, we can prove this result for $\rho (y)dy$ almost
every $y$. So $(u^{i}(s,X_{s}^{t,y}),(\sigma ^{*}\nabla
u^{i})(s,X_{s}^{t,y}))$ solves the BSDE$(g(X_{T}^{t,y}),f)$, i.e. $\rho
(y)dy $ a.s., we have
\[
u^{i}(s,X_{s}^{t,y})=g(X_{T}^{t,y})+%
\int_{s}^{T}f(s,X_{s}^{t,y},u^{i}(s,X_{s}^{t,y}),\sigma ^{*}\nabla
u^{i}(s,X_{s}^{t,y}))ds-\int_{s}^{T}(\sigma ^{*}\nabla
u^{i})(r,X_{r}^{t,y})dB_{r}.
\]
Then by the uniqueness of the BSDE, we know $%
u^{1}(s,X_{s}^{t,y})=u^{2}(s,X_{s}^{t,x})$ and $(\sigma ^{*}\nabla
u^{1})(s,X_{s}^{t,y})=(\sigma ^{*}\nabla u^{2})(s,X_{s}^{t,y})$.
Taking $s=t$ we deduce that $u^{1}(t,y)=u^{2}(t,y)$, $ dt\otimes
dy$-a.s. $\hfill \Box$

\section{Sobolev's solution for PDE with obstacle under
monotonicity condition}
\sectionmark{Sobolev Solutions for PDE with obstacle}

In this section we study the PDE with obstacle associated with
$(g,f,h)$, which satisfy the assumptions
\ref{term-pde}-\ref{bar-pde} for $n=1 $. We will prove the existence
and uniqueness of a weak solution to the obstacle problem. We will
restrict our study to the case when $\varphi $ is polynomial
increasing in $y$, i.e.
\begin{assumption}\label{poly} We assume that for some $\kappa _{1}\in %
\mathbb{R}$, $\beta _{1}>0,\forall y\in \mathbb{R}$,
\[
\left| \varphi (y)\right| \leq \kappa _{1}(1+\left| y\right|
^{\beta _{1}}).
\]
\end{assumption}
For the sake of PDE with obstacle, we introduce the reflected BSDE
associated with $(g,f,h)$, like in El Karoui et al. \cite{EKPPQ}:
\begin{equation}
\label{rbsde1}
 \left\lbrace
\begin{split}
& Y_{s}^{t,x}
 =g(X_{T}^{t,x})+
\int_{s}^{T}f(r,X_{r}^{t,x},Y_{r}^{t,x},Z_{r}^{t,x})dr+K_{T}^{t,x}-K_{t}^{t,x}-\int_{s}^{T}Z_{s}^{t,x}dB_{s},\;
P\text{-}a.s \; \forall \,  s \in [t,T]  \\
& Y_{s}^{t,x} \geq L_{s}^{t,x}, \; P\text{-}a.s \\
& \int_{t}^{T}(Y_{s}^{t,x}-L_{s}^{t,x})dK_{s}^{t,x}=0, \, \, P\text{-}a.s.\\
\end{split}
\right.
\end{equation}
where $L_{s}^{t,x}=h(s,X_{s}^{t,x})$ is a continuous process.
Moreover following Lepeltier et al \cite{LMX},  we shall need to
estimate
\begin{eqnarray*}
E[\sup_{t\leq s\leq T}\varphi ^{2}(e^{\mu t}(L_{s}^{t,x})^{+})]
&=&E[\sup_{t\leq s\leq T}\varphi ^{2}(e^{\mu t}h(s,X_{s}^{t,x})^{+})] \\
&\leq &Ce^{2\beta _{1}\mu T}E[\sup_{t\leq s\leq T}(1+\left|
X_{s}^{t,x}\right| ^{2\beta _{1}\beta })] \\
&\leq &C(1+\left| x\right| ^{2\beta _{1}\beta }),
\end{eqnarray*}
where $C$ is a constant which can be changed line by line. By
assumption \ref{bar-pde}-(ii), with same techniques we get for $x\in %
\mathbb{R}$, $E[\sup_{t\leq s\leq T}\varphi
^{2}((L_{s}^{t,x})^{+})]<+\infty $. Thanks to the assumption
\ref{term-pde} and \ref{coef-pde}, by the equivalence of norms
\ref{equiv1} and \ref{equiv2}, we have
\[
g(X_{T}^{t,x}) \in \mathbf{L}^2({\cal F}_T) \mbox{ and }
f(s,X_{s}^{t,x},0,0) \in \mathbf{H}^2(0,T).
\]
 By the
existence and uniqueness theorem for the RBSDE in \cite{LMX}, for each $%
(t,x)$, there exists a unique triple $(Y^{t,x},Z^{t,x},K^{t,x})\in \mathbf{S}%
^{2}(t,T)\times \mathbf{H}_{d}^{2}(t,T)\times \mathbf{A}^{2}(t,T)$
of $\{\mathcal{F}_{s}^{t}\}$ progressively measurable processes,
which is the solution of the reflected
BSDE with parameters $(g(X_{T}^{t,x}),f(s,X_{s}^{t,x},y,z),$ $h(s,X_{s}^{t,x}))$%
We shall give the probabilistic interpretation for the solution of
PDE with obstacle (\ref{OPDE}).

The main result of this section is

\begin{theorem}
\label{mr2}Assume that assumptions \ref{term-pde}-\ref{mark} hold
and $\rho (x)=(1+\left| x\right| )^{-p}$ with $p\geq \gamma $ where
$\gamma =\beta _{1}\beta +\beta +d+1$. There exists a pair $(u,\nu
)$, which is the solution of the PDE with obstacle (\ref{OPDE})
associated to $(g,f,h)$ i.e. $(u,\nu )$ satisfies
Definition \ref{o-pde}-(i) -(iii). Moreover the solution is given by: $%
u(t,x)=Y_{t}^{t,x}$, a.e. where
$(Y_{s}^{t,x},Z_{s}^{t,x},K_{s}^{t,x})_{t\leq s\leq T}$ is the
solution of RBSDE (\ref{rbsde1}), and
\begin{equation}
Y_{s}^{t,x}=u(s,X_{s}^{t,x}),Z_{s}^{t,x}=(\sigma ^{*}\nabla
u)(s,X_{s}^{t,x}).  \label{con-pre}
\end{equation}
Moreover, we have for every measurable bounded and positive
functions $\phi $ and $\psi $,
\begin{equation}
\int_{\mathbb{R}^{d}}\int_{t}^{T}\phi (s,\widehat{X}_{s}^{t,x})J(\widehat{X}%
_{s}^{t,x})\psi (s,x)1_{\{u=h\}}(s,x)d\nu (s,x)=\int_{\mathbb{R}%
^{d}}\int_{t}^{T}\phi (s,x)\psi (s,X_{s}^{t,x})dK_{s}^{t,x}\text{, a.s..}
\label{con-k}
\end{equation}
If $(\overline{u},\overline{\nu })$ is another solution of the PDE (\ref
{OPDE}) such that $\overline{\nu }$ satisfies (\ref{con-k}) with some $%
\overline{K}$ instead of $K$, where $\overline{K}$ is a continuous process
in $\mathbf{A}_{\mathcal{F}}^{2}(t,T)$, then $\overline{u}=u$ and $\overline{%
\nu }=\nu $.
\end{theorem}

\begin{remark}
The expression (\ref{con-k}) gives us the probabilistic
interpretation (Feymamn-Kac's formula) for the measure $\nu $ via
the increasing process $K^{t,x}$ of the RBSDE. This formula was
first introduced in Bally et al. \cite{BCEF}, where the authors
prove (\ref{con-k}) when $f$ is Lipschitz on $y$ and $z$ uniformly
in $(t,\omega )$. Here we generalize their result to the case when
$f$ is monotonic in $y$ and Lipschitz in $z$.
\end{remark}

\textbf{Proof.} As in the proof of theorem \ref{mr1} in section 4, we first
notice that $(u,\nu )$ solves (\ref{OPDE}) if and only if
\[
(\widehat{u}(t,x),d\widehat{\nu }(t,x))=(e^{\mu t}u(t,x),e^{\mu t}d\nu (t,x))
\]
is the solution of the PDE with obstacle $(\widehat{g},\widehat{f},\widehat{h%
})$, where $\widehat{g}$, $\widehat{f}$ are defined as in
(\ref{tran-f}) with
\[
\widehat{h}(t,x)=e^{\mu t}h(t,x).
\]
Then the coefficient $\widehat{f}$ satisfies the same assumptions in
assumption \ref{coef-pde1}\textbf{\ }with (iii) replaced by (\ref{tran-f}),
which means that $f$ is decreasing on $y$ in the $1$-dimensional case. The
obstacle $\widehat{h}$ still satisfies assumption \ref{bar-pde}, for $\mu =0$%
. In the following we will use $(g,f,h)$ instead of $(\widehat{g},\widehat{f}%
,\widehat{h})$, and suppose that $(g,f,h)$ satisfies assumption \ref
{term-pde}, \ref{coef-pde}, \ref{bar-pde}, \ref{mark} and \ref{coef-pde1}%
\textbf{\ }with (iii) replaced by (\ref{tran-f}).\\[0.2cm]
\textbf{ \large a) Existence} : \textit{The existence of a solution
will be proved in 4 steps. From step 1 to step 3, we suppose that
$f$ does not depend on $\nabla u,$ satisfies
assumption \ref{coef-pde1}' for $n=1$, and $f(t,x,0)\in \mathbf{L}%
^{2}([0,T]\times \mathbb{R}^{d},dt\otimes \rho (x)dx)$\textbf{.} In the step
4, we study the case when $f$ depend on $\nabla u$}.\\[0.2cm]
\textbf{\large Step 1} :  \textit{Suppose $g(x)$, $f(t,x,0)$,
$h^{+}(t,x)$ uniformly bounded} i.e. that there exists a constant
$C$ such that
\[
\left| g(x)\right| +\sup_{0\leq t\leq T}\left| f(t,x,0)\right| +\sup_{0\leq
t\leq T}h^{+}(t,x)\leq C.
\]

We will use the penalization method. For $n\in \mathbb{N}$, we consider for
all $s\in [t,T]$,
\[
Y_{s}^{n,t,x}=g(X_{T}^{t,x})+\int_{s}^{T}f(r,X_{r}^{t,x},Y_{r}^{n,t,x})dr+n%
\int_{s}^{T}(Y_{r}^{n,t,x}-h(r,X_{r}^{t,x}))^{-}dr-%
\int_{s}^{T}Z_{r}^{n,t,x}dB_{r}.
\]
From Theorem \ref{mr1} in section 3, we know that $u_{n}(t,x):=Y_{t}^{n,t,x}$%
, is solution of the PDE$(g,f_{n})$, where $%
f_{n}(t,x,y,x)=f(t,x,y,z)+n(y-h(t,x))^{-}$, i.e. for every $\phi \in
C_{c}^{1,\infty }([0,T]\times \mathbb{R}^{d})$%
\begin{equation}
\label{o-equa1} \begin{split}
 \int_{t}^{T}(u_{s}^{n},\partial
_{s}\phi )ds & +(u^{n}(t,\cdot ),\phi
(t,\cdot ))-(g(\cdot ),\phi (\cdot ,T))+\int_{t}^{T}\mathcal{E}%
(u_{s}^{n},\phi _{s})ds  \nonumber \\
&=\int_{t}^{T}(f(s,\cdot ,u_{s}^{n}),\phi
_{s})ds+n\int_{t}^{T}((u^{n}-h)^{-}(s,\cdot ),\phi _{s})ds.
\end{split}
\end{equation}
Moreover
\begin{equation}
Y_{s}^{n,t,x}=u_{n}(s,X_{s}^{t,x}),Z_{s}^{n,t,x}=\sigma ^{*}\nabla
u_{n}(s,X_{s}^{t,x}),  \label{rep1}
\end{equation}
Set $K_{s}^{n,t,x}=n \displaystyle
\int_{t}^{s}(Y_{r}^{n,t,x}-h(r,X_{r}^{t,x}))^{-}dr$. Then by
(\ref{rep1}), we have that $K_{s}^{n,t,x}=n \displaystyle
\int_{t}^{s}(u_{n}-h)^{-}(r,X_{r}^{t,x})dr$.

Following the estimates and convergence results for
$(Y^{n,t,x},Z^{n,t,x})$ in the step 1 of the proof of Theorem 2.2 in
\cite{LMX}, for $m$, $n\in \mathbb{N}$, we have,as $m,n\rightarrow
\infty $
\[
E\int_{t}^{T}\left| Y_{s}^{n,t,x}-Y_{s}^{m,t,x}\right|
^{2}ds+E\int_{t}^{T}\left| Z_{s}^{n,t,x}-Z_{s}^{m,t,x}\right|
^{2}ds+E\sup_{t\leq s\leq T}\left|
K_{s}^{n,t,x}-K_{s}^{m,t,x}\right| ^{2}\rightarrow 0,
\]
and
\begin{equation*}
\sup_n E\int_{0}^{T}(\left| Y_{s}^{n,t,x}\right| ^{2}+\left|
Z_{s}^{n,t,x}\right|^2+ (K_{T}^{n,t,x})^{2})\leq C.
\end{equation*}
By the equivalence of the norms (\ref{equiv2}), we get
\begin{eqnarray*}
&&\int_{\mathbb{R}^{d}}\int_{t}^{T}\rho (x)(\left|
u_{n}(s,x)-u_{m}(s,x)\right| ^{2}+\left| \sigma ^{*}\nabla u_{n}(s,x)-\sigma
^{*}\nabla u_{m}(s,x)\right| ^{2})dsdx \\
&\leq &\frac{1}{k_{2}}\int_{\mathbb{R}^{d}}\rho (x)E\int_{t}^{T}(\left|
Y_{s}^{n,t,x}-Y_{s}^{m,t,x}\right| ^{2}+\left|
Z_{s}^{n,t,x}-Z_{s}^{m,t,x}\right| ^{2})dsdx\rightarrow 0.
\end{eqnarray*}
Thus $(u_{n})$ is a Cauchy sequence in $\mathcal{H}$, and the limit $%
u=\lim_{n\rightarrow \infty }u_{n}$ belongs to  $\mathcal{H}$.

Denote $\nu _{n}(dt,dx)=n(u_{n}-h)^{-}(t,x)dtdx$ and $\pi _{n}(dt,dx)=\rho
(x)\nu _{n}(dt,dx)$, then by (\ref{equiv1})
\begin{eqnarray*}
\pi _{n}([0,T]\times \mathbb{R}^{d}) &=&\int_{\mathbb{R}^{d}}\int_{0}^{T}%
\rho (x)\nu _{n}(dt,dx)=\int_{\mathbb{R}^{d}}\int_{0}^{T}\rho
(x)n(u_{n}-h)^{-}(t,x)dtdx \\
&\leq &\frac{1}{k_{2}}\int_{\mathbb{R}^{d}}\rho (x)E\left|
K_{T}^{n,0,x}\right| dx\leq C\int_{\mathbb{R}^{d}}\rho (x)dx<\infty .
\end{eqnarray*}
It follows that
\begin{equation}
\sup_{n}\pi _{n}([0,T]\times \mathbb{R}^{d})<\infty .  \label{est-measure}
\end{equation}
In the same way like in the existence proof step 2 of theorem 14
in \cite {BCEF}, we can prove that $\pi _{n}([0,T]\times
\mathbb{R}^{d})$ is bounded and then $\pi _{n}$ is tight. So we
may pass to a subsequence and get $\pi _{n}\rightarrow \pi $ where
$\pi $ is a positive measure. Define $\nu =\rho
^{-1}\pi $; $\nu $ is a positive measure such that $\int_{0}^{T}\int_{%
\mathbb{R}^{d}}\rho (x)d\nu (t,x)<\infty $, and so we have for
$\phi \in C_{c}^{1,\infty}([0,T]\times \mathbb{R}^{d})$ with
compact support in $x$,
\[
\int_{\mathbb{R}^{d}}\int_{t}^{T}\phi d\nu _{n}=\int_{\mathbb{R}%
^{d}}\int_{t}^{T}\frac{\phi }{\rho }d\pi _{n}\rightarrow \int_{\mathbb{R}%
^{d}}\int_{t}^{T}\frac{\phi }{\rho }d\pi =\int_{\mathbb{R}%
^{d}}\int_{t}^{T}\phi d\nu .
\]

Now passing to the limit in the PDE$(g,f_{n})$, we check that $(u,\nu )$
satisfies the PDE with obstacle $(g,f,h)$, i.e. for every $\phi \in
C_{c}^{1,\infty }([0,T]\times \mathbb{R}^{d})$, we have
\begin{eqnarray}
&&\int_{t}^{T}(u_{s},\partial _{s}\phi )ds+(u(t,\cdot ),\phi (t,\cdot
))-(g(\cdot ),\phi (\cdot ,T))+\int_{t}^{T}\mathcal{E}(u_{s},\phi _{s})ds
\nonumber \\
&=&\int_{t}^{T}(f(s,\cdot ,u_{s}),\phi _{s})ds+\int_{t}^{T}\int_{\mathbb{R}%
^{d}}\phi (s,x)1_{\{u=h\}}(s,x)d\nu (x,s).  \label{equa1}
\end{eqnarray}

The last is to prove that $\nu $ satisfies the probabilistic interpretation (%
\ref{con-k}). Since $K^{n,t,x}$ converges to $K^{t,x}$ uniformly in $t$, the
measure $dK^{n,t,x}\rightarrow dK^{t,x}$ weakly in probability.

Fix two continuous functions $\phi $, $\psi $ : $[0,T]\times \mathbb{R}%
^{d}\rightarrow \mathbb{R}^{+}$ which have compact support in $x$ and a
continuous function with compact support $\theta :\mathbb{R}^{d}\rightarrow %
\mathbb{R}^{+}$, we have
\begin{eqnarray*}
&&\int_{\mathbb{R}^{d}}\int_{t}^{T}\phi (s,\widehat{X}_{s}^{t,x})J(\widehat{X%
}_{s}^{t,x})\psi (s,x)\theta (x)d\nu (s,x) \\
&=&\lim_{n\rightarrow \infty }\int_{\mathbb{R}^{d}}\int_{t}^{T}\phi (s,%
\widehat{X}_{s}^{t,x})J(\widehat{X}_{s}^{t,x})\psi (s,x)\theta
(x)n(u_{n}-h)^{-}(t,x)dtdx \\
&=&\lim_{n\rightarrow \infty }\int_{\mathbb{R}^{d}}\int_{t}^{T}\phi
(s,x)\psi (s,X_{s}^{t,x})\theta
(X_{s}^{t,x})n(u_{n}-h)^{-}(t,X_{s}^{t,x})dtdx \\
&=&\lim_{n\rightarrow \infty }\int_{\mathbb{R}^{d}}\int_{t}^{T}\phi
(s,x)\psi (s,X_{s}^{t,x})\theta (X_{s}^{t,x})dK_{s}^{n,t,x}dx \\
&=&\int_{\mathbb{R}^{d}}\int_{t}^{T}\phi (s,x)\psi (s,X_{s}^{t,x})\theta
(X_{s}^{t,x})dK_{s}^{t,x}dx.
\end{eqnarray*}

We take $\theta =\theta _{R}$ to be the regularization of the indicator
function of the ball of radius $R$ and pass to the limit with $R\rightarrow
\infty $, it follows that
\begin{equation}\label{con-k1}
\int_{\mathbb{R}^{d}}\int_{t}^{T}\phi (s,\widehat{X}_{s}^{t,x})J(\widehat{X}%
_{s}^{t,x})\psi (s,x)d\nu (s,x)=\int_{\mathbb{R}^{d}}\int_{t}^{T}\phi
(s,x)\psi (s,X_{s}^{t,x})dK_{s}^{t,x}dx.
\end{equation}
Since $(Y_{s}^{n,t,x},Z_{s}^{n,t,x},K_{s}^{n,t,x})$ converges to $%
(Y_{s}^{t,x},Z_{s}^{t,x},K_{s}^{t,x})$ as $n\rightarrow \infty $ in $\mathbf{%
S}^{2}(t,T)\mathbf{\times H}^{2}(t,T)\times \mathbf{A}^{2}(t,T)$, and $%
(Y_{s}^{t,x},Z_{s}^{t,x},K_{s}^{t,x})$ is the solution of RBSDE$%
(g(X_{T}^{t,x}),f,h)$, then we have
\[
\int_{t}^{T}(Y_{s}^{t,x}-L_{s}^{t,x})dK_{s}^{t,x}=%
\int_{t}^{T}(u-h)(t,X_{s}^{t,x})dK_{s}^{t,x}=0,\text{a.s.}
\]
it follows that $dK_{s}^{t,x}=1_{\{u=h\}}(s,X_{s}^{t,x})dK_{s}^{t,x}$. In (%
\ref{con-k1}), setting $\psi =1_{\{u=h\}}$ yields
\[
\int_{\mathbb{R}^{d}}\int_{t}^{T}\phi (s,\widehat{X}_{s}^{t,x})J(\widehat{X}%
_{s}^{t,x})1_{\{u=h\}}(s,x)d\nu (s,x)=\int_{\mathbb{R}^{d}}\int_{t}^{T}\phi
(s,\widehat{X}_{s}^{t,x})J(\widehat{X}_{s}^{t,x})d\nu (s,x)\text{, a.s.}
\]
Note that the family of functions $A(\omega )=\{(s,x)\rightarrow \phi (s,%
\widehat{X}_{s}^{t,x}):\phi \in C_{c}^{\infty }\}$ is an algebra which
separates the points (because $x\rightarrow \widehat{X}_{s}^{t,x}$ is a
bijection). Given a compact set $G$, $A(\omega )$ is dense in $C([0,T]\times
G)$. It follows that $J(\widehat{X}_{s}^{t,x})1_{\{u=h\}}(s,x)d\nu (s,x)=J(%
\widehat{X}_{s}^{t,x})d\nu (s,x)$ for almost every $\omega $. While $J(%
\widehat{X}_{s}^{t,x})>0$ for almost every $\omega $, we get $d\nu
(s,x)=1_{\{u=h\}}(s,x)d\nu (s,x)$, and (\ref{con-k}) follows.

Then we get easily that $Y_{s}^{t,x}=u(s,X_{s}^{t,x})$ and $%
Z_{s}^{t,x}=\sigma ^{*}\nabla u(s,X_{s}^{t,x})$, in view of the convergence
results for $(Y_{s}^{n,t,x},Z_{s}^{n,t,x})$ and the equivalence of the
norms. So $u(s,X_{s}^{t,x})=Y_{s}^{t,x}\geq h(t,x)$. Specially for $s=t$, we
have $u(t,x)\geq h(t,x)$\\[0.2cm]
\textbf{\large Step 2} :  \textit{As in the proof of the RBSDE in
Theorem 2.2 in \cite{LMX}, step 2, we relax the bounded condition on
the barrier $h$ in step 1, and prove the existence of the solution
under assumption \ref{bar-pde}}.

Similarly to step 2 in the proof of theorem 2.2 in \cite{LMX},
after some transformation, we know that it is sufficient to prove
the existence of the solution for the PDE with obstacle $(g,f,h)$,
where $(g,f,h)$ satisfies
\[
g(x),f(t,x,0)\leq 0.
\]

Let $h(t,x)$ satisfy assumption \ref{bar-pde} for $\mu =0$, i.e. $\forall
(t,x)\in [0,T]\times ,\mathbb{R}^{d}$
\[
\varphi (h(t,x)^{+})\in \mathbf{L}^{2}(\mathbb{R}^{d};\rho (x)dx),
\]
and
\[
\left| h(t,x)\right| \leq \kappa (1+\left| x\right| ^{\beta }).
\]
Set
\[
h_{n}(t,x)=h(t,x)\wedge n,
\]
then the function $h_{n}(t,x)$ are continuous, $\sup_{0\leq t\leq
T}h_{n}^{+}(t,x)\leq n$, and $h_{n}(s,X_{s}^{t,x})\rightarrow
h(s,X_{s}^{t,x})$ in $\mathbf{S}_{\mathcal{F}}^{2}(t,T)$, in view of Dini's
theorem and dominated convergence theorem.

We consider the PDE with obstacle associated with $(g,f,h_{n})$. By the
results of step 1, there exists $(u_{n},\nu _{n})$, which is the solution of
the PDE with obstacle associated to $(g,f,h_{n})$, where $u_{n}\in \mathcal{H%
}$ and $\nu _{n}$ is a positive measure such that $\int_{0}^{T}\int_{%
\mathbb{R}^{d}}\rho (x)d\nu _{n}(t,x)<\infty $. Moreover
\begin{equation}\label{boun-bar1}
\left\{
\begin{split}
& Y_{s}^{n,t,x}=u_{n}(s,X_{s}^{t,x}),Z_{s}^{n,t,x}=\sigma ^{*}\nabla
u_{n}(s,X_{s}^{t,x}), \\
& \int_{\mathbb{R}^{d}}\int_{t}^{T}\phi (s,\widehat{X}_{s}^{t,x})J(\widehat{X%
}_{s}^{t,x})\psi (s,x)1_{\{u_{n}=h_{n}\}}(s,x)d\nu _{n}(s,x)=\int_{\mathbb{R}%
^{d}}\int_{t}^{T}\phi (s,x)\psi (s,X_{s}^{t,x})dK_{s}^{n,t,x}dx,
\end{split}
\right.
\end{equation}
Here $(Y^{n,t,x},Z^{n,t,x},K^{n,t,x})$ is the solution of the RBSDE$%
(g(X_{T}^{t,x}),f,h_{n})$. Thanks to proposition \ref{esti} in
Appendix, and the bounded assumption of $g$ and $f$, we know that
\begin{equation}
\label{est-kb}
\begin{split}
&  E \, \big[\int_{t}^{T}\big(\left| Y_{s}^{n,t,x}\right| ^{2}+
\left| Z_{s}^{n,t,x}\right| ^{2}\, \big)ds+(K_{T}^{n,0,x})^{2}
\big]\\ & \leq C \, \big(1+E \, \big[\varphi ^{2}(\sup_{0\leq t\leq
T}h^{+}(t,X_{t}^{0,x})\, \big)+\sup_{0\leq t\leq
T}(h^{+}(t,X_{t}^{0,x}))^{2}\big] \big) \\
 &\leq C(1+\left| x\right|
^{2\beta _{1}\beta }+\left| x\right| ^{2\beta }).\\
\end{split}
\end{equation}
By the Lemma 2.3 in \cite{LMX}, $Y_{s}^{n,t,x}%
\rightarrow Y_{s}^{t,x}$ in $\mathbf{S}^{2}(0,T)$, $%
Z_{s}^{n,t,x}\rightarrow Z_{s}^{t,x}$ in $\mathbf{H}_{d}^{2}(0,T)$
and $K_{s}^{n,t,x}\rightarrow K_{s}^{t,x}$ in $\mathbf{A}^{2}(0,T)$, as $n\rightarrow \infty $.
Moreover $%
(Y_{s}^{t,x},Z_{s}^{t,x},K_{s}^{t,x})$ is the solution of RBSDE$%
(g(X_{T}^{t,x}),f,h)$.

By the convergence result of $(Y_{s}^{n,t,x},Z_{s}^{n,t,x})$ and
the equivalence of the norms (\ref{equiv2}), we get
\begin{eqnarray*}
&&\int_{\mathbb{R}^{d}}\rho (x)\int_{t}^{T}(\left|
u_{n}(t,x)-u_{m}(t,x)\right| ^{2}+\left| \sigma ^{*}\nabla u_{n}(s,x)-\sigma
^{*}\nabla u_{m}(s,x)\right| ^{2})dsdx \\
&\leq &\frac{1}{k_{2}}\int_{\mathbb{R}^{d}}\rho (x)E\int_{t}^{T}(\left|
Y_{s}^{n,t,x}-Y_{s}^{m,t,x}\right| ^{2}+\left|
Z_{s}^{n,t,x}-Z_{s}^{m,t,x}\right| ^{2})dsdx\rightarrow 0.
\end{eqnarray*}
So $\{u_{n}\}$ is a Cauchy sequence in $\mathcal{H}$, and admits a limit $%
u\in \mathcal{H}$. Moreover
$Y_{s}^{t,x}=u(s,X_{s}^{t,x}),Z_{s}^{t,x}=\sigma ^{*}\nabla
u(s,X_{s}^{t,x})$. In particular $u(t,x)=Y_{t}^{t,x}\geq h(t,x)$.

Set $\pi _{n}=\rho \nu _{n}$, like in step 1, we first need to prove that $%
\pi _{n}([0,T]\times \mathbb{R}^{d})$ is uniformly bounded. In (\ref
{boun-bar1}), let $\phi =\rho $, $\psi =1$, then we have
\[
\int_{\mathbb{R}^{d}}\int_{0}^{T}\rho (\widehat{X}_{s}^{0,x})J(\widehat{X}%
_{s}^{0,x})d\nu _{n}(s,x)=\int_{\mathbb{R}^{d}}\int_{0}^{T}\rho
(x)dK_{s}^{n,0,x}dx.
\]

Recall Lemma \ref{flo-bon}: there exist two constants $c_{1}>0$ and $c_{2}>0$
such that $\forall x\in \mathbb{R}^{d}$, $0\leq t\leq T$%
\[
c_{1}\leq E\left( \frac{\rho (t,\widehat{X}_{t}^{0,x})J(\widehat{X}%
_{t}^{0,x})}{\rho (x)}\right) \leq c_{2}.
\]
Applying H\"{o}lder's inequality and Schwartz's inequality, we have
\begin{eqnarray*}
&&\pi _{n}([0,T]\times \mathbb{R}^{d}) \\
&=&\int_{\mathbb{R}^{d}}\int_{0}^{T}\rho (x)\nu _{n}(dt,dx) \\
&=&\int_{\mathbb{R}^{d}}\int_{0}^{T}\frac{\rho ^{\frac{1}{2}}(x)}{\rho ^{%
\frac{1}{2}}(t,\widehat{X}_{t}^{0,x})J^{\frac{1}{2}}(\widehat{X}_{t}^{0,x})}%
\rho ^{\frac{1}{2}}(x)\rho ^{\frac{1}{2}}(t,\widehat{X}_{t}^{0,x})J^{\frac{1%
}{2}}(\widehat{X}_{t}^{0,x})\nu _{n}(dt,dx) \\
&\leq &E[\left( \int_{\mathbb{R}^{d}}\int_{0}^{T}\frac{\rho (x)}{\rho (t,%
\widehat{X}_{t}^{0,x})J(\widehat{X}_{t}^{0,x})}\rho (x)\nu
_{n}(dt,dx)\right) ^{\frac{1}{2}}\left( \int_{\mathbb{R}^{d}}\int_{0}^{T}%
\rho (t,\widehat{X}_{t}^{0,x})J(\widehat{X}_{t}^{0,x})\nu _{n}(dt,dx)\right)
^{\frac{1}{2}}] \\
&\leq &\left( E\int_{\mathbb{R}^{d}}\int_{0}^{T}\frac{\rho (x)}{\rho (t,%
\widehat{X}_{t}^{0,x})J(\widehat{X}_{t}^{0,x})}\rho (x)\nu
_{n}(dt,dx)\right) ^{\frac{1}{2}}\left( E\int_{\mathbb{R}^{d}}\int_{0}^{T}%
\rho (t,\widehat{X}_{t}^{0,x})J(\widehat{X}_{t}^{0,x})\nu _{n}(dt,dx)\right)
^{\frac{1}{2}} \\
&=&\left( \int_{\mathbb{R}^{d}}\int_{0}^{T}E\left( \frac{\rho (x)}{\rho (t,%
\widehat{X}_{t}^{0,x})J(\widehat{X}_{t}^{0,x})}\right) \rho (x)\nu
_{n}(dt,dx)\right) ^{\frac{1}{2}}\left( \int_{\mathbb{R}^{d}}E%
\int_{0}^{T}dK_{t}^{n,0,x}\rho (x)dx\right) ^{\frac{1}{2}} \\
&\leq &\left( \frac{1}{c_{1}}\int_{\mathbb{R}^{d}}\int_{0}^{T}\rho (x)\nu
_{n}(dt,dx)\right) ^{\frac{1}{2}}\left( \int_{\mathbb{R}^{d}}\rho
(x)E[K_{T}^{n,0,x}]dx\right) ^{\frac{1}{2}}.
\end{eqnarray*}
So by (\ref{est-kb}) and (\ref{equiv1}), we get
\begin{eqnarray}
\sup_{n}\pi _{n}([0,T]\times \mathbb{R}^{d}) &\leq &C\int_{\mathbb{R}%
^{d}}\rho (x)E[K_{T}^{n,0,x}]dx  \label{est-measure2} \\
&\leq &C\int_{\mathbb{R}^{d}}\rho (x)(1+\left| x\right| ^{\beta _{1}\beta
}+\left| x\right| ^{\beta })dx<\infty .  \nonumber
\end{eqnarray}
Using the same arguments as in step 1, we deduce that $\pi _{n}$ is tight.
So we may pass to a subsequence and get $\pi _{n}\rightarrow \pi $ where $%
\pi $ is a positive measure.

Define $\nu =\rho ^{-1}\pi $, then $\nu $ is a positive measure such that $%
\int_{0}^{T}\int_{\mathbb{R}^{d}}\rho (x)d\nu (t,x)<\infty $. Then for $\phi
\in C([0,T]\times \mathbb{R}^{d})$ with compact support in $x$, we have as $%
n\rightarrow \infty $,
\[
\int_{t}^{T}\int \phi d\nu _{n}=\int_{t}^{T}\int \frac{\phi }{\rho }d\pi
_{n}\rightarrow \int_{t}^{T}\int \frac{\phi }{\rho }d\pi =\int_{t}^{T}\int
\phi d\nu .
\]

Now passing to the limit in the PDE$(g,f,h_{n})$, we check that $(u,\nu )$
satisfies the PDE with obstacle associated to $(g,f,h)$, i.e. for every $%
\phi \in C_{c}^{1,\infty }([0,T]\times \mathbb{R}^{d})$%
\begin{eqnarray}
&&\int_{t}^{T}(u_{s},\partial _{s}\phi )ds+(u(t,\cdot ),\phi (t,\cdot
))-(g(\cdot ),\phi (\cdot ,T))+\int_{t}^{T}\mathcal{E}(u_{s},\phi _{s})ds
\nonumber \\
&=&\int_{t}^{T}(f(s,\cdot ,u_{s}),\phi _{s})ds+\int_{t}^{T}\int_{\mathbb{R}%
^{d}}\phi (s,x)1_{\{u=h\}}d\nu (x,s).  \label{equa2}
\end{eqnarray}

Then we will check if the probabilistic interpretation (\ref{con-k}) still
holds. Fix two continuous functions $\phi $, $\psi $ : $[0,T]\times %
\mathbb{R}^{d}\rightarrow \mathbb{R}^{+}$ which have compact support in $x$.
With the convergence result of $K^{n,t,x}$, which implies $%
dK^{n,t,x}\rightarrow $ $dK^{t,x}$ weakly in probability, in the same way as
step 1, passing to the limit in (\ref{boun-bar1}) we have
\begin{equation}
\int_{\mathbb{R}^{d}}\int_{t}^{T}\phi (s,\widehat{X}_{s}^{t,x})J(\widehat{X}%
_{s}^{t,x})\psi (s,x)d\nu (s,x)=\int_{\mathbb{R}^{d}}\int_{t}^{T}\phi
(s,x)\psi (s,X_{s}^{t,x})dK_{s}^{t,x}dx  \nonumber  \label{con-k2}
\end{equation}
Since $(Y_{s}^{t,x},Z_{s}^{t,x},K_{s}^{t,x})$ is the solution of RBSDE$%
(g(X_{T}^{t,x}),f,h)$, then by the integral condition, we deduce the $%
dK_{s}^{t,x}=1_{\{u=h\}}(s,X_{s}^{t,x})dK_{s}^{t,x}$. In (\ref{con-k2}),
setting $\psi =1_{\{u=h\}}$ yields
\[
\int_{\mathbb{R}^{d}}\int_{t}^{T}\phi (s,\widehat{X}_{s}^{t,x})J(\widehat{X}%
_{s}^{t,x})1_{\{u=h\}}(s,x)d\nu (s,x)=\int_{\mathbb{R}^{d}}\int_{t}^{T}\phi
(s,\widehat{X}_{s}^{t,x})J(\widehat{X}_{s}^{t,x})d\nu (s,x).
\]
With the same arguments, we get that $d\nu (s,x)=1_{\{u=h\}}(s,x)d\nu (s,x)$%
, and (\ref{con-k})holds for $\nu $ and $K$.\\[0.2cm]
\textbf{\large Step 3} :  \textit{Now we will relax the bounded condition on $g(x)$ and $%
f(t,x,0)$}. Then for $m,n\in \mathbb{N}$, let
\begin{eqnarray*}
g_{m,n}(x) &=&(g(x)\wedge n)\vee (-m), \\
f_{m,n}(t,x,y) &=&f(t,x,y)-f(t,x,0)+(f(t,x,0)\wedge n)\vee (-m).
\end{eqnarray*}
\textbf{\ }So $g_{m,n}(x)$ and $f_{m,n}(t,x,0)$ are bounded and for fixed $%
m\in \mathbb{N}$, as $n\rightarrow \infty $, we have
\begin{eqnarray*}
g_{m,n}(x) &\rightarrow &g_{m}(x)\ \text{in }\mathbf{L}^{2}\mathbb{(R}%
^{d},\rho (x)dx), \\
f_{m,n}(t,x,0) &\rightarrow &f_{m}(t,x,0)\text{ in }\mathbf{L}^{2}\mathbf{(}%
[0,T]\mathbb{\times R}^{d},dt\otimes \rho (x)dx),
\end{eqnarray*}
where
\begin{eqnarray*}
g_{m}(x) &=&g(x)\vee (-m), \\
f_{m}(t,x,y) &=&f(t,x,y)-f(t,x,0)+f(t,x,0)\vee (-m).
\end{eqnarray*}
Then as $m\rightarrow \infty $, we have
\begin{eqnarray*}
g_{m}(x) &\rightarrow &g(x)\ \text{in
}\mathbf{L}^{2}\mathbb{(R}^{d},\rho
(x)dx), \\
f_{m}(t,x,0) &\rightarrow &f(t,x,0)\text{ in }\mathbf{L}^{2}\mathbf{(}[0,T]%
\mathbb{\times R}^{d},dt\otimes \rho (x)dx),
\end{eqnarray*}
in view of assumption \ref{term-pde} and $f(t,x,0)\in \mathbf{L}%
^{2}([0,T]\times \mathbb{R}^{d},dt\otimes \rho (x)dx).$

Now we consider the PDE with obstacle associated to
$(g_{m,n},f_{m,n},h)$. By step 2, there exists a $(u_{m,n},\nu
_{m,n})$ which is the solution of the PDE with obstacle associated
to $(g_{m,n},f_{m,n},h)$. In particular the representation
 formulas (\ref{con-pre}) and (\ref{con-k}) are satisfied.  Denote by $(Y^{m,n,t,x},Z^{m,n,t,x},K^{m,n,t,x})$
 the solution of the RBSDE $(g_{m,n}(X_{T}^{t,x}),f_{m,n},h)$.

Recall the convergence results in step 3 of theorem 2.2 in
\cite{LMX}, we know
that for fixed $m\in \mathbb{N}$, as $n\rightarrow \infty $, $%
(Y_{s}^{m,n,t,x},Z_{s}^{m,n,t,x},K_{s}^{m,n,t,x})\rightarrow
(Y_{s}^{m,t,x},Z_{s}^{m,t,x},K_{s}^{m,t,x})$ in
$\mathbf{S}^{2}(0,T)\times \mathbf{H}_{d}^{2}(0,T)\times
\mathbf{A}^{2}(0,T)$, and that $%
(Y_{s}^{m,t,x},Z_{s}^{m,t,x},K_{s}^{m,t,x})$ is the solution of RBSDE$%
(g_{m}(X_{T}^{t,x}),f_{m},h)$.

By It\^{o}'s formula, we have for $n,p\in \mathbb{N}$,
\begin{eqnarray*}
&&E\int_{t}^{T}(\left| Y_{s}^{m,n,t,x}-Y_{s}^{m,p,t,x}\right| ^{2}+\left|
Z_{s}^{m,n,t,x}-Z_{s}^{m,p,t,x}\right| ^{2})ds \\
&\leq &CE\left| g_{m,n}(X_{T}^{t,x})-g_{m,p}(X_{T}^{t,x})\right|
^{2}+CE\int_{t}^{T}\left|
f_{m,n}(s,X_{s}^{t,x},0)-f_{m,p}(s,X_{s}^{t,x},0)\right| ^{2}ds,
\end{eqnarray*}
so by the equivalence of the norms (\ref{equiv1}) and (\ref{equiv2}), it
follows that as $n\rightarrow \infty $,
\begin{eqnarray*}
&&\int_{\mathbb{R}^{d}}\int_{t}^{T}\rho (x)(\left|
u_{m,n}(t,x)-u_{m,p}(t,x)\right| ^{2}+\left| \sigma ^{*}\nabla
u_{m,n}(s,x)-\sigma ^{*}\nabla u_{m,p}(s,x)\right| ^{2})dsdx \\
&\leq &\frac{Ck_{1}}{k_{2}}\int_{\mathbb{R}^{d}}\rho (x)\left|
g_{m,n}(x)-g_{m,p}(x)\right| ^{2}dx+\frac{Ck_{1}}{k_{2}}\int_{\mathbb{R}%
^{d}}\int_{t}^{T}\rho (x)\left| f_{m,n}(s,x,0)-f_{m,p}(s,x,0)\right|
^{2}dsdx\rightarrow 0.
\end{eqnarray*}
i.e. for each fixed $m\in \mathbb{N}$, $\{u_{m,n}\}$ is a Cauchy sequence in
$\mathcal{H}$, and admits a limit $u_{m}\in \mathcal{H}$. Moreover $%
Y_{s}^{m,t,x}=u_{m}(s,X_{s}^{t,x}),Z_{s}^{m,t,x}=\sigma ^{*}\nabla
u_{m}(s,X_{s}^{t,x})$, a.s., in particular $u_{m}(t,x)=Y_{t}^{m,t,x}\geq
h(t,x)$.

Then we find the measure $\nu _{m}$ by the sequence $\{\nu _{m,n}\}$. Set $%
\pi _{m,n}=\rho \nu _{m,n}$, by proposition \ref{esti} in
Appendix, we have for each $m,n\in \mathbb{N}$, $0\leq t\leq T$
\begin{eqnarray}
\ \ \ E(\left| K_{T}^{m,n,t,x}\right| ^{2}) &\leq
&CE[g_{m,n}^{2}(X_{T}^{t,x})+\int_{0}^{T}f_{m,n}^{2}(s,X_{s}^{t,x},0,0)ds+%
\varphi ^{2}(\sup_{t\leq s\leq T}(h^{+}(s,X_{s}^{t,x})))  \nonumber \\
&&+\sup_{t\leq s\leq T}(h^{+}(s,X_{s}^{t,x}))^{2}+1+\varphi ^{2}(2T)]
\nonumber \\
&\leq &CE[g(X_{T}^{t,x})^{2}+\int_{0}^{T}f^{2}(s,X_{s}^{t,x},0,0)ds+\varphi
^{2}(\sup_{0\leq s\leq T}(h^{+}(s,X_{s}^{t,x})))  \nonumber \\
&&+\sup_{0\leq s\leq T}(h^{+}(s,X_{s}^{t,x}))^{2}+1+\varphi ^{2}(2T)]
\nonumber \\
&\leq &C(1+\left| x\right| ^{2\beta _{1}\beta }+\left| x\right| ^{2\beta }).
\label{est-kmn}
\end{eqnarray}
By the same way as in step 2, we deduce that for each fixed $m\in \mathbb{N}$%
, $\pi _{m,n}$ is tight, we may pass to a subsequence and get $\pi
_{m,n}\rightarrow \pi _{m}$ where $\pi _{m}$ is a positive measure. If we
define $\nu _{m}=\rho ^{-1}\pi _{m}$, then $\nu _{m}$ is a positive measure
such that $\int_{0}^{T}\int_{\mathbb{R}^{d}}\rho (x)d\nu _{m}(t,x)<\infty $.
So we have for all $\phi \in C([0,T]\times \mathbb{R}^{d})$ with compact
support in $x$,
\[
\int_{t}^{T}\int \phi d\nu _{m,n}=\int_{t}^{T}\int \frac{\phi }{\rho }d\pi
_{m,n}\rightarrow \int_{t}^{T}\int \frac{\phi }{\rho }d\pi
_{m}=\int_{t}^{T}\int \phi d\nu _{m}.
\]

Now for each fixed $m\in \mathbb{N}$, let $n\rightarrow \infty $, in the PDE$%
(g_{m,n},f_{m,n},h)$, we check that $(u_{m},\nu _{m})$ satisfies the PDE
with obstacle associated to $(g_{m},f_{m},h)$, and by the weak convergence
result of $dK^{m,n,t,x}$, we have easily that the probabilistic
interpretation (\ref{con-k}) holds for $\nu _{m}$ and $K^{m}$.

Then let $m\rightarrow \infty $, by the convergence results in step 4 of
theorem 2.2 in \cite{LMX}, we apply the same method as before. We deduce that $%
\lim_{m\rightarrow \infty }u_{m}=u$ is in $\mathcal{H}$ and $%
Y_{s}^{t,x}=u(s,X_{s}^{t,x}),Z_{s}^{t,x}=\sigma ^{*}\nabla u(s,X_{s}^{t,x})$%
, a.s., where $(Y^{t,x},Z^{t,x},K^{t,x})$ is the solution of the RBSDE$%
(g,f,h)$, in particular, setting $s=t$, $u_{(}t,x)=Y_{t}^{t,x}\geq h(t,x)$.

From (\ref{est-kmn}), it follows that
\[
E[(K_{T}^{m,t,x})^{2}]\leq C(1+\left| x\right| ^{2\beta \beta _{1}}+\left|
x\right| ^{2\beta }).
\]
By the same arguments, we can find the measure $\nu $ by the sequence $\{\nu
_{m}\}$, which satisfies that for all $\phi $ and $\psi $ with compact
support,
\[
\int_{\mathbb{R}^{d}}\int_{t}^{T}\phi (s,\widehat{X}_{s}^{t,x})J(\widehat{X}%
_{s}^{t,x})\psi (s,x)1_{\{u=h\}}(s,x)d\nu (s,x)=\int_{\mathbb{R}%
^{d}}\int_{t}^{T}\phi (s,x)\psi (s,X_{s}^{t,x})dK_{s}^{t,x}dx.
\]

Finally we find a solution $(u,\nu)$ to the PDE with obstacle $(g,f,h)$,
when $f$ does not depend on $\nabla u$. So for every $\phi \in
C_{c}^{1,\infty }([0,T]\times \mathbb{R}^{d})$%
\begin{eqnarray}
&&\int_{t}^{T}(u_{s},\partial _{s}\phi )ds+(u(t,\cdot ),\phi (t,\cdot
))-(g(\cdot ),\phi (\cdot ,T))+\int_{t}^{T}\mathcal{E}(u_{s},\phi _{s})ds
\nonumber \\
&=&\int_{t}^{T}(f(s,\cdot ,u_{s}),\phi _{s})ds+\int_{t}^{T}\int_{\mathbb{R}%
^{d}}\phi (s,x)1_{\{u=h\}}d\nu (x,s).  \label{equa4}
\end{eqnarray}

\textbf{\large Step 4 } : \textit{Finally we study the case when $f$
depends on $\nabla u$, and satisfies a Lipschitz condition on
$\nabla u$}. We construct a mapping $\Psi $ from $\mathcal{H}$ into itself. For some $%
\overline{u}\in \mathcal{H}$, define
\[
u=\Psi (\overline{u}),
\]
where $(u,\nu )$ is a weak solution of the PDE with obstacle $%
(g,f(t,x,u,\sigma \nabla \overline{u}),h)$. Then by this mapping, we denote
a sequence $\{u_{n}\}$ in $\mathcal{H}$, beginning with a function $v^{0}\in
\mathbf{L}^{2}([0,T]\times \mathbb{R}^{d},dt\otimes \rho (x)dx)$. Since $%
f(t,x,u,v^{0}(t,x))$ satisfies the assumptions of step 3, the PDE$%
(g,f(t,x,u,v^{0}(t,x)),h)$ admits a solution $(u_{1},v_{1})\in \mathcal{H}$.
For $n\in \mathbb{N}$, set $u_{n}(t,x)=\Psi (u_{n-1}(t,x))$.

Symmetrically we introduce a mapping $\Phi $ from
$\mathbf{H}^{2}(t,T)\times \mathbf{H}_{d}^{2}(t,T)$ into
itself. For $V^{t,x,0}=v^{0}(s,X_{s}^{t,x}))$, then $V_{s}^{t,x}\in \mathbf{H%
}_{d}^{2}(t,T)$ in view of the equivalence of the norms. Set
\[
(Y^{t,x,n},Z^{t,x,n})=\Phi (Y^{t,x,n-1},Z^{t,x,n-1}),
\]
where $(Y^{t,x,n},Z^{t,x,n},K^{t,x,n})$is the solution of the RBSDE with
parameters $g(X_{T}^{t,x})$, $f(s,X_{s}^{t,x},Y_{s}^{t,x},Z_{s}^{t,x,n-1})$
and $h(s,X_{s}^{t,x})$.Then $Y_{s}^{t,x,n}=u_{n}(s,X_{s}^{t,x})$, $%
Z_{s}^{t,x,n}=\sigma ^{*}\nabla u_{n}(s,X_{s}^{t,x})$ a.s. and
\[
\int_{\mathbb{R}^{d}}\int_{t}^{T}\phi (s,\widehat{X}_{s}^{t,x})J(\widehat{X}%
_{s}^{t,x})\psi (s,x)1_{\{u=h\}}(s,x)d\nu _{n}(s,x)=\int_{\mathbb{R}%
^{d}}\int_{t}^{T}\phi (s,x)\psi (s,X_{s}^{t,x})dK_{s}^{t,x,n}dx.
\]

Set $\widetilde{u}_{n}(t,x):=u_{n}(t,x)-u_{n-1}(t,x)$. To deal
with the difference $\widetilde{u}_{n}$, we need the difference of
the corresponding
BSDE, denote $\widetilde{Y}_{s}^{t,x,n}:=Y_{s}^{t,x,n}-Y_{s}^{t,x,n-1}$, $%
\widetilde{Z}_{s}^{t,x,n}:=Z_{s}^{t,x,n}-Z_{s}^{t,x,n-1}$,
$\widetilde{K}_{s}^{t,x,n}:=K_{s}^{t,x,n}-K_{s}^{t,x,n-1}$. It
follows from It\^{o}'s formula, for some $\alpha$, $\gamma \in
\mathbb{R}$,
\begin{eqnarray*}
\ e^{\gamma t}E\left| \widetilde{Y}_{s}^{t,x,n}\right| ^{2}
&+&E\int_{s}^{T}e^{\gamma r}(\gamma \left|
\widetilde{Y}_{r}^{t,x,n}\right|
^{2}+\left| \widetilde{ Z}_{r}^{t,x,n}\right| ^{2})dr\  \\
&\leq &E\int_{s}^{T}e^{\gamma r}(\frac{k^{2}}{\alpha }\left|
\widetilde{ Y}_{r}^{t,x,n}\right| ^{2}+\alpha \left|
\widetilde{Z}_{r}^{t,x,n-1}\right| ^{2})dr,
\end{eqnarray*}
since
\begin{eqnarray*}
&&\int_{s}^{T}e^{\gamma r}\widetilde{ Y}_{r}^{t,x,n}d\widetilde{ K}_{r}^{t,x,n} \\
&=&\int_{s}^{T}e^{\gamma
r}(Y_{s}^{t,x,n}-h(r,X_{r}^{t,x}))dK^{t,x,n}+\int_{s}^{T}e^{\gamma
r}(Y_{s}^{t,x,n-1}-h(r,X_{r}^{t,x}))dK^{t,x,n-1} \\
&&\ -\int_{s}^{T}e^{\gamma
r}(Y_{s}^{t,x,n}-h(r,X_{r}^{t,x}))dK^{t,x,n-1}+\int_{s}^{T}e^{\gamma
r}(Y_{s}^{t,x,n-1}-h(r,X_{r}^{t,x}))dK^{t,x,n} \\
\  &\leq &0.
\end{eqnarray*}
then by the equivalence of the norms, for $\gamma =1+\frac{2k_{1}^{2}}{%
k_{2}^{2}}k^{2}$, we have
\begin{eqnarray*}
&&\ \int_{\mathbb{R}^{d}}\int_{t}^{T}e^{\gamma s}(\left|
\widetilde{ u}_{n}(s,x)\right| ^{2}+\left| \sigma ^{*}\nabla
(\widetilde{u}_{n})(s,x)\right| ^{2})\rho (x)dsdx \\
\  &\leq
&(\frac{1}{2})^{n-1}\int_{\mathbb{R}^{d}}\int_{t}^{T}e^{\gamma
s}(\left| \widetilde{ u}_{2}(s,x)\right| ^{2}+\left| \sigma
^{*}\nabla
(\widetilde{ u}_{2})(s,x)\right| ^{2})\rho (x)dsdx \\
\  &\leq &(\frac{1}{2})^{n-1}(\left\| u_{1}(s,x)\right\| _{\gamma
}^{2}+\left\| u_{2}(s,x)\right\| _{\gamma }^{2}).
\end{eqnarray*}
where $\left\| u\right\| _{\gamma }^{2}:=\int_{\mathbb{R}^{d}}%
\int_{t}^{T}e^{\gamma s}(\left| u(s,x)\right| ^{2}+\left| \sigma ^{*}\nabla
u(s,x)\right| ^{2})\rho (x)dsdx$, which is equivalent to the norm $\left\|
\cdot \right\| $ of $\mathcal{H}$. So $\{u_{n}\}$ is a Cauchy sequence in $%
\mathcal{H}$, it admits a limit $u$ in $\mathcal{H}$, which is the solution
to the PDE with obstacle (\ref{Pde1}). Then consider $\sigma ^{*}\nabla u$
as a known function by the result of step 3, we know that there exists a
positive measure $\nu $ such that $\int_{0}^{T}\int_{\mathbb{R}^{d}}\rho
(x)d\nu (t,x)<\infty $, and for every $\phi \in C_{c}^{1,\infty
}([0,T]\times \mathbb{R}^{d})$,
\begin{eqnarray}
&&\int_{t}^{T}(u_{s},\partial _{s}\phi )ds+(u(t,\cdot ),\phi (t,\cdot
))-(g(\cdot ),\phi (\cdot ,T))+\int_{t}^{T}\mathcal{E}(u_{s},\phi _{s})ds
\nonumber \\
&=&\int_{t}^{T}(f(s,\cdot ,u_{s},\sigma ^{*}\nabla u_{s}),\phi
_{s})ds+\int_{t}^{T}\int_{\mathbb{R}^{d}}\phi (s,x)1_{\{u=h\}}d\nu (x,s).
\label{equa5}
\end{eqnarray}
Moreover, for $t\leq s\leq T$,
\[
Y_{s}^{t,x}=u(s,X_{s}^{t,x}),Z_{s}^{t,x}=\sigma ^{*}\nabla u(s,X_{s}^{t,x}),%
\text{ a.s.a.e.,}
\]
and
\begin{eqnarray*}
&&\int_{\mathbb{R}^{d}}\int_{t}^{T}\phi (s,\widehat{X}_{s}^{t,x})J(\widehat{X%
}_{s}^{t,x})\psi (s,x)1_{\{u=h\}}(s,x)d\nu (s,x) \\
&=&\int_{\mathbb{R}^{d}}\int_{t}^{T}\phi (s,x)\psi
(s,X_{s}^{t,x})dK_{s}^{t,x}.
\end{eqnarray*}
\\[0.2cm]

\textbf{ \large b) Uniqueness } :  Set $(\overline{u},\overline{\nu
})$ to be another solution of the PDE with
obstacle (\ref{OPDE}) associated to $(g,f,h)$; with $\overline{\nu }$ verifies (%
\ref{con-k}) for an increasing process $\overline{K}$. We fix $\phi :%
\mathbb{R}^{d}\rightarrow \mathbb{R}$, a smooth function in $C_{c}^{2}(%
\mathbb{R}^{d})$ with compact support and denote $\phi _{t}(s,x)=\phi (%
\widehat{X}_{s}^{t,x})J(\widehat{X}_{s}^{t,x})$. From proposition \ref{test1}%
, one may use $\phi _{t}(s,x)$ as a test function in the PDE$(g,f,h)$ with $%
\partial _{s}\phi (s,x)ds$ replaced by a stochastic integral with respect to
the semimartingale $\phi _{t}(s,x)$. Then we get, for $t\leq s\leq T$%
\begin{eqnarray}
&&\int_{\mathbb{R}^{d}}\int_{s}^{T}\overline{u}(r,x)d\phi _{t}(r,x)dx+(%
\overline{u}(s,\cdot ),\phi _{t}(s,\cdot ))-(g(\cdot ),\phi _{t}(\cdot
,T))+\int_{s}^{T}\mathcal{E}(\overline{u}_{r},\phi _{r})dr  \label{o-pde-u1}
\\
&=&\int_{s}^{T}\int_{\mathbb{R}^{d}}f(r,x,\overline{u}(r,x),\sigma
^{*}\nabla \overline{u}(r,x))\phi _{t}(r,\cdot )dr+\int_{s}^{T}\int_{\mathbf{%
R}^{d}}\phi _{t}(r,x)1_{\{\overline{u}=h\}}d\overline{\nu }(x,r).  \nonumber
\end{eqnarray}
By (\ref{dphi}) in Lemma \ref{comp}, we have
\begin{eqnarray*}
\int_{\mathbb{R}^{d}}\int_{s}^{T}\overline{u}d_{r}\phi _{t}(r,x)dx
&=&\int_{s}^{T}(\int_{\mathbb{R}^{d}}(\sigma ^{*}\nabla \overline{u}%
)(r,x)\phi _{t}(r,x)dx)dB_{r} \\
&&+\int_{s}^{T}\int_{\mathbb{R}^{d}}\left( (\sigma ^{*}\nabla \overline{u}%
)(\sigma ^{*}\nabla \phi _{r})+\phi _{t}\nabla ((\frac{1}{2}\sigma
^{*}\nabla \sigma +b)\overline{u})\right) dxdr.
\end{eqnarray*}
Substitute this equality in (\ref{o-pde-u1}), we get
\begin{eqnarray*}
\int_{\mathbb{R}^{d}}\overline{u}(s,x)\phi _{t}(s,x)dx &=&(g(\cdot ),\phi
_{t}(\cdot ,T))-\int_{s}^{T}(\int_{\mathbb{R}^{d}}(\sigma ^{*}\nabla
\overline{u})(r,x)\phi _{t}(r,x)dx)dB_{r} \\
&&+\int_{\mathbb{R}^{d}}\int_{s}^{T}f(r,x,\overline{u}(r,x),\sigma
^{*}\nabla \overline{u}(r,x))\phi _{t}(s,\cdot )dr+\int_{s}^{T}\int_{\mathbf{%
R}^{d}}\phi _{t}(r,x)1_{\{\overline{u}=h\}}d\overline{\nu }(x,r).
\end{eqnarray*}
Then by changing of variable $y=\widehat{X}_{r}^{t,x}$ and applying (\ref
{con-k}) for $\overline{\nu }$, we obtain
\begin{eqnarray*}
&&\int_{\mathbb{R}^{d}}\overline{u}(s,X_{s}^{t,y})\phi (y)dy \\
&=&\int_{\mathbb{R}^{d}}g(X_{T}^{t,y})\phi (y)dy+\int_{s}^{T}\phi
(y)f(s,X_{s}^{t,y},\overline{u}(s,X_{s}^{t,y}),\sigma ^{*}\nabla \overline{u}%
(s,X_{s}^{t,y})ds \\
&&+\int_{s}^{T}\int_{\mathbb{R}^{d}}\phi (y)1_{\{\overline{u}%
=h\}}(r,X_{s}^{t,y})d\overline{K}_{r}^{t,y}dy-\int_{s}^{T}(\int_{\mathbb{R}%
^{d}}(\sigma ^{*}\nabla \overline{u})(r,X_{r}^{t,y})\phi (y)dy)dB_{r}.
\end{eqnarray*}
Since $\phi $ is arbitrary, we can prove that for $\rho (y)dy$ almost every $%
y$, $(\overline{u}(s,X_{s}^{t,y}),(\sigma ^{*}\nabla \overline{u}%
)(s,X_{s}^{t,y}),\widehat{K}_{s}^{t,x})$ solves the RBSDE$%
(g(X_{T}^{t,y}),f,h)$. Here $\widehat{K}_{s}^{t,x}=\int_{t}^{s}1_{\{%
\overline{u}=h\}}(r,X_{r}^{t,y})d\overline{K}_{r}^{t,y}$. Then by the
uniqueness of the solution of the reflected BSDE, we know $\overline{u}%
(s,X_{s}^{t,y})=Y_{s}^{t,y}=u(s,X_{s}^{t,x})$, $(\sigma ^{*}\nabla \overline{%
u})(s,X_{s}^{t,y})=Z_{s}^{t,y}=(\sigma ^{*}\nabla u)(s,X_{s}^{t,y})$ and $%
\widehat{K}_{s}^{t,y}=K_{s}^{t,y}$. Taking $s=t$ we deduce that $\overline{u}%
(t,y)=u(t,y)$, $\rho (y)dy$-a.s. and by the probabilistic interpretation (%
\ref{con-k}), we obtain
\[
\int_{s}^{T}\int \phi _{t}(r,x)1_{\{\overline{u}=h\}}(r,x)d\overline{\nu }%
(x,r)=\int_{s}^{T}\int \phi _{t}(r,x)1_{\{u=h\}}(r,x)d\nu (x,r).
\]
So $1_{\{\overline{u}=h\}}(r,x)d\overline{\nu
}(x,r)=1_{\{u=h\}}(r,x)d\nu (x,r)$.$\hfill \Box$
\section{Appendix}
\sectionmark{Appendix}
\subsection{Proof of proposition \ref{test1}}

 First we consider the case when
$f$ does not depend on $z$ and satisfies
assumption \ref{coef-pde1}\textbf{'}. As in step 2 of the proof of theorem \ref{mr1}, we approximate $g$ and $%
f $ as in (\ref{appr}), then $g_{n}\rightarrow g$ in $\mathbf{L}^{2}(%
\mathbb{R}^{d},\rho (x)dx)$ and $f_{n}(t,x,0)\rightarrow f(t,x,0)$ in $%
\mathbf{L}^{2}([0,T]\times \mathbb{R}^{d},dt\otimes \rho (x)dx)$, as $%
n\rightarrow \infty $.

Since for each $n\in \mathbb{N}$, $\left| g_{n}\right| \leq n$ and $\left|
f_{n}(t,x,0)\right| \leq n$, by the result of the step 1 of theorem \ref{mr1}%
, the PDE$(g_{n},f_{n})$ admits the weak solution $u_{n}\in \mathcal{H}$ and
$\sup_{0\leq t\leq T}\left| u_{n}(t,x)\right| \leq C_{n}$. So we know
\[
\left| f_{n}(t,x,u_{n}(t,x))\right| ^{2}\leq \left| f_{n}(t,x,0)\right|
^{2}+\varphi (\sup_{0\leq t\leq T}\left| u_{n}(t,x)\right| )\leq C_{n}.
\]
Set $F_{n}(t,x):=f_{n}(t,x,u_{n}(t,x))$, then $F_{n}(t,x)\in \mathbf{L}%
^{2}([0,T]\times \mathbb{R}^{n},dt\otimes \rho (x)dx).$

From proposition 2.3 in Bally and Matoussi \cite{BM2001}, for
$\phi \in C_{c}^{2}(\mathbb{R}^{d})$, we get, for $t\leq s\leq T$%
\begin{equation}
\label{tesf-b}
\begin{split}
&\int_{\mathbb{R}^{d}}\int_{s}^{T}u_{n}(r,x)d\phi
_{t}(r,x)dx+(u_{n}(s,\cdot ),\phi _{t}(s,\cdot ))-(g_{n}(\cdot
),\phi _{t}(\cdot ,T))+\int_{s}^{T}\mathcal{E}(u_{n}(r,\cdot ),\phi
_{t}(r,\cdot
))dr  \\
&=\int_{\mathbb{R}^{d}}\int_{s}^{T}f(r,x,u_{n}(r,x))\phi
_{t}(r,x)drdx+\int_{\mathbb{R}^{d}}\int_{s}^{T}(f_{n}(r,x,0)-f(r,x,0))\phi
_{t}(r,x)drdx.
\end{split}
\end{equation}
By step 2, we know that as $n\rightarrow \infty $, $u_{n}\rightarrow
u$ in $\mathcal{H}$, where $u$ is a weak solution of the PDE$(g,f)$,
i.e.
\begin{eqnarray*}
u_{n} &\rightarrow &u\text{ in }\mathbf{L}^{2}([0,T]\times \mathbb{R}%
^{d},dt\otimes \rho (x)dx), \\
\sigma ^{*}\nabla u_{n} &\rightarrow &\nabla u\text{ in }\mathbf{L}%
^{2}([0,T]\times \mathbb{R}^{d},dt\otimes \rho (x)dx).
\end{eqnarray*}
Then there exists a function $u^{*}$ in $\mathbf{L}^{2}([0,T]\times %
\mathbb{R}^{d},dt\otimes \rho (x)dx)$, such that for a subsequence of $%
\{u_{n}\}$, $\left| u_{n_{k}}\right| \leq \left| u^{*}\right| $ and $%
u_{n_{k}}\rightarrow u$, $dt\otimes dx$-a.e. Thanks to assumption
\ref {coef-pde1}\textbf{'}-(iii), we have that $
f(r,x,u_{n}(r,x))\rightarrow f(r,x,u(r,x))$, $dt\otimes dx$-a.e.
Now, for all compact support function $\phi \in
C_{c}^{2}(\mathbb{R}^{d})$,  the second term in the right hand side
of (\ref{tesf-b}) converge to 0 as $ n \to \infty$ and it is not
hard to prove by using the dominated convergence theorem the term in
the left hand side of (\ref{tesf-b}) converges. Thus, we conclude
that $\lim_{n\rightarrow
\infty}\int_{\mathbb{R}^{d}}\int_{s}^{T}f(r,x,u_{n}(r,x))\phi
_{t}(r,x)drdx $ exists. Moreover by the monotonocity condition of
$f$  and the same arguments as in step 2 of the proof of theorem
\ref{mr1}, we get for all compact support function $\phi \in
C_{c}^{2}(\mathbb{R}^{d})$
\begin{eqnarray*}
&&\int_{\mathbb{R}^{d}}\int_{s}^{T}u(r,x)d\phi _{t}(r,x)dx+(u(s,\cdot ),\phi
_{t}(s,\cdot ))-(g(\cdot ),\phi _{t}(\cdot ,T))+\int_{s}^{T}\mathcal{E}%
(u(r,\cdot ),\phi _{t}(r,\cdot ))dr \\
&=& \int_{\mathbb{R}^{d}}\int_{s}^{T}f(r,x,u(r,x))\phi _{t}(r,x)drdx
\, .
\end{eqnarray*}

Now we consider the case when $f$ depends on $\nabla u$ and
satisfies the assumption \ref{coef-pde1}\textbf{\ }with (iii)
replaced by (\ref{tran-f}). Like in the step 3 of the proof of
theorem \ref{mr1}, we construct a mapping $\Psi $ from $\mathcal{H}$
into
itself. Then by this mapping, we define a sequence $\{u_{n}\}$ in $\mathcal{H%
}$, beginning with a matrix-valued function $v^{0}\in \mathbf{L}%
^{2}([0,T]\times \mathbb{R}^{n\times d},dt\otimes \rho (x)dx)$. Since $%
f(t,x,u,v^{0}(t,x))$ satisfies the assumptions of step 2, the PDE$%
(g,f(t,x,u,v^{0}(t,x)))$ admits a unique solution $u_{1}\in \mathcal{H}$.
For $n\in \mathbb{N}$, denote
\[
u_{n}(t,x)=\Psi (u_{n-1}(t,x)),
\]
i.e. $u_{n}$ is the weak solution of the PDE$(g,f(t,x,u,\sigma
^{*}\nabla u_{n-1}(t,x)))$. Set
$\widetilde{u}_{n}(t,x):=u_{n}(t,x)-u_{n-1}(t,x)$. In
order to estimate the difference, we introduce the corresponding BSDE$%
(g,f_{n})$ for $n=1$, where $f_{n}(t,x,u)=f(t,x,u,\nabla
u_{n-1}(t,x))$. So we have
$Y_{s}^{n,t,x}=u_{n}(s,X_{s}^{t,x}),Z_{s}^{n,t,x}=\sigma \nabla
u_{n}(s,X_{s}^{t,x})$. Then we apply the It\^{o}'s formula to
$|\widetilde{ Y}^{n,t,x}|^{2}$, where
$\widetilde{Y}_{s}^{n,t,x}:=Y_{s}^{n,t,x}-Y_{s}^{n-1,t,x}
$. With the equivalence of the norms, similarly as in step 3, for $\gamma =1+%
\frac{2k_{1}^{2}}{k_{2}^{2}}k^{2}$, we have
\begin{eqnarray*}
&&\int_{\mathbb{R}^{d}}\int_{t}^{T}e^{\gamma s}(\left|
\widetilde{u}_{n}(s,x)\right| ^{2}+\left| \sigma ^{*}\nabla
(\widetilde{u}_{n})(s,x)\right| ^{2})\rho (x)dsdx \\
&\leq
&(\frac{1}{2})^{n-1}\int_{\mathbb{R}^{d}}\int_{t}^{T}e^{\gamma
s}(\left| \widetilde{ u}_{2}(s,x)\right| ^{2}+\left| \sigma
^{*}\nabla
(\widetilde{ u}_{2})(s,x)\right| ^{2})\rho (x)dsdx \\
&\leq &(\frac{1}{2})^{n-1}(\left\| u_{1}(s,x)\right\| _{\gamma }^{2}+\left\|
u_{2}(s,x)\right\| _{\gamma }^{2}).
\end{eqnarray*}
where $\left\| u\right\| _{\gamma }^{2}:=\int_{\mathbb{R}^{d}}%
\int_{t}^{T}e^{\gamma s}(\left| u(s,x)\right| ^{2}+\left| \sigma ^{*}\nabla
u(s,x)\right| ^{2})\rho (x)dsdx$, which is equivalent to the norm $\left\|
\cdot \right\| $ of $\mathcal{H}$. So $\{u_{n}\}$ is a Cauchy sequence in $%
\mathcal{H}$, it admits a limit $u$ in $\mathcal{H}$, and by the fixed point
theorem, $u$ is a solution of the PDE$(g,f)$.

Then for each $n\in \mathbb{N}$, we have for $\phi \in C_{c}^{2}(\mathbb{R}%
^{d})$%
\begin{eqnarray*}
&&\int_{\mathbb{R}^{d}}\int_{s}^{T}u_{n}(r,x)d\phi
_{t}(r,x)dx+(u_{n}(s,\cdot ),\phi _{t}(s,\cdot ))-(g(\cdot ),\phi _{t}(\cdot
,T))+\int_{s}^{T}\mathcal{E}(u_{n}(r,\cdot ),\phi _{r}(r,\cdot ))dr \\
&=&\int_{\mathbb{R}^{d}}\int_{s}^{T}f(r,x,u_{n}(r,x),\sigma ^{*}\nabla
u_{n-1}(r,x))\phi _{t}(r,x)drdx \\
&=&\int_{\mathbb{R}^{d}}\int_{s}^{T}f(r,x,u_{n}(r,x),\sigma ^{*}\nabla
u(r,x))\phi _{t}(r,x)drdx \\
&&+\int_{\mathbb{R}^{d}}\int_{s}^{T}[f(r,x,u_{n}(r,x),\sigma ^{*}\nabla
u_{n-1}(r,x))-f(r,x,u_{n}(r,x),\sigma ^{*}\nabla u(r,x))]\phi _{t}(r,x)drdx.
\end{eqnarray*}
Noticing that $f$ is Lipschitz in $z$, we get
\[
\left| f(r,x,u_{n}(r,x),\sigma ^{*}\nabla
u_{n-1}(r,x))-f(r,x,u_{n}(r,x),\sigma ^{*}\nabla u(r,x))\right| \leq k\left|
\sigma ^{*}\nabla u_{n-1}(r,x)-\sigma ^{*}\nabla u(r,x)\right| .
\]
So the last term of the right side converges to $0$, since $\{\sigma
^{*}\nabla u_{n}\}$ converges to $\sigma ^{*}\nabla u$ in $\mathbf{L}%
^{2}([0,T]\times \mathbb{R}^{d},dt\otimes \rho (x)dx)$. Now we are in the
same situation as in the first part of proof, and in the same way, we deduce
that the following holds: for $\phi \in C_{c}^{2}(\mathbb{R}^{d})$%
\begin{eqnarray*}
&&\int_{\mathbb{R}^{d}}\int_{s}^{T}u(r,x)d\phi _{t}(r,x)dx+(u(s,\cdot ),\phi
_{t}(s,\cdot ))-(g(\cdot ),\phi _{t}(\cdot ,T))+\int_{s}^{T}\mathcal{E}%
(u(r,\cdot ),\phi _{t}(r,\cdot ))dr \\
&=&\int_{\mathbb{R}^{d}}\int_{s}^{T}f(r,x,u(r,x),\sigma ^{*}\nabla
u(r,x))\phi _{t}(r,x)drdx,dt\otimes dx\text{, a.s..}
\end{eqnarray*}

Now if $f$ satisfies assumption \ref{coef-pde1}, we know  that $u$ is solution of the PDE$(g,f)$ if and only if $%
\widehat{u}=e^{\mu t}u$ is solution of the PDE$(\widehat{g},\widehat{f})$,
where
\[
\widehat{g}(x)=e^{\mu T}g(x),\;\widehat{f}(t,x,y,x)=e^{\mu t}f(t,x,e^{-\mu
t}y,e^{-\mu t}z)-\mu y,
\]
and $\widehat{f}$ satisfies assumption \ref{coef-pde1}-(iii)
replaced by (\ref{tran-f}). So we know now: for $\phi \in
C_{c}^{2}(\mathbb{R}^{d})$,
\begin{eqnarray*}
&&\int_{\mathbb{R}^{d}}\int_{s}^{T}\widehat{u}(r,x)d\phi _{t}(r,x)dx+(%
\widehat{u}(s,\cdot ),\phi _{t}(s,\cdot ))-(\widehat{g}(\cdot ),\phi
_{t}(\cdot ,T))+\int_{s}^{T}\mathcal{E}(\widehat{u}(r,\cdot ),\phi
_{t}(r,\cdot ))dr \\
&=&\int_{\mathbb{R}^{d}}\int_{s}^{T}\widehat{f}(r,x,\widehat{u}(r,x),\nabla
\widehat{u}(r,x))\phi _{t}(r,x)drdx,dt\otimes dx\text{, a.s..}
\end{eqnarray*}
Notice that $ d(e^{\mu r}u(r,x))=\mu e^{\mu r}u(r,x)dr+e^{\mu
r}d(u(r,x))$, so  by the integration by parts formula (for
stochastic process), we get
\begin{eqnarray*}
&&\int_{\mathbb{R}^{d}}\int_{s}^{T}u(r,x)d\phi _{t}(r,x)dx \\
&=&\int_{\mathbb{R}^{d}}\int_{s}^{T}e^{-\mu r}\widehat{u}(r,x)d\phi
_{t}(r,x)dxdr \\
&=&e^{-\mu T}(\widehat{g}(\cdot ),\phi _{t}(\cdot ,T))-e^{-\mu s}(\widehat{u}%
(s,\cdot ),\phi _{t}(s,\cdot ))+\mu \int_{s}^{T}e^{-\mu r}\int_{\mathbb{R}%
^{d}}\widehat{u}(r,x)\phi _{t}(r,x)dxdr \\
&&-\int_{s}^{T}\int_{\mathbb{R}^{d}}e^{-\mu r}\phi _{t}(r,x)[\mathcal{L}%
\widehat{u}(r,x)+\widehat{f}(r,x,\widehat{u}(r,x),\nabla \widehat{u}%
(r,x))]drdx.
\end{eqnarray*}
Using (\ref{tran-m}), we get that for $\phi \in C_{c}^{2}(\mathbb{R}^{d})$,
\begin{eqnarray*}
\int_{\mathbb{R}^{d}}\int_{s}^{T}u(r,x)d\phi _{t}(r,x)dx &=&(g(\cdot ),\phi
_{t}(\cdot ,T))-(u(s,\cdot ),\phi _{t}(s,\cdot ))-\int_{s}^{T}\int_{%
\mathbb{R}^{d}}\phi _{t}(r,x)\mathcal{L}u(r,x)drdx \\
&&+\int_{s}^{T}\int_{\mathbb{R}^{d}}\phi _{t}(r,x)f(r,x,u(r,x),\nabla
u(r,x))drdx \\
&=&(g(\cdot ),\phi _{t}(\cdot ,T))-(u(s,\cdot ),\phi _{t}(s,\cdot
))-\int_{s}^{T}\mathcal{E(}u(r,\cdot ),\phi _{t}(r,\cdot ))dr \\
&&+\int_{s}^{T}\int_{\mathbb{R}^{d}}\phi _{t}(r,x)f(r,x,u(r,x),\nabla
u(r,x))drdx,
\end{eqnarray*}
and finally, the result follows. $\hfill \Box$

\subsection{Some a priori estimates}

In this subsection, we consider the non-markovian Reflected BSDE
associated to $(\xi ,f,L)$ :
\begin{equation*}
\left\lbrace
\begin{split}
&Y_{t} =\xi
+\int_{t}^{T}f(t,Y_{s},Z_{s})ds+K_{T}-K_{t}-\int_{t}^{T}Z_{s}dB_{s}, \\
&Y_{t} \geq L_{t},\\
& \int_{0}^{T}(Y_{s}-L_{s})dK_{s}=0 \\
\end{split}
\right.
\end{equation*}
under the following assumptions :\\[0.2cm]
\textbf{(H1)} a final condition $\xi \in
\mathbf{L}^{2}(\mathcal{F}_{T})$,\\[0.1cm]
\textbf{(H2)} a coefficient $f:\Omega \times [0,T]\times
\mathbb{R}\times \mathbb{R}^{d}\rightarrow \mathbb{R} $, which is
such that for some continuous increasing function $\varphi
:\mathbb{R}_{+}\longrightarrow \mathbb{R}_{+}$, a real numbers $\mu
$ and $C>0$:\\[0.1cm]
$
\begin{array}{cl}
\text{(i)} & f(\cdot ,y,z)\text{ is progressively measurable,
}\forall
(y,z)\in \mathbb{R}\times \mathbb{R}^{d}; \\
\text{(ii)} & \left| f(t,y,0)\right| \leq \left| f(t,0,0)\right|
+\varphi (\left| y\right| )\text{, }\forall (t,y)\in [0,T]\times
\mathbb{R}\text{,
a.s.;} \\
\text{(iii)} & E\int_{0}^{T}\left| f(t,0,0)\right| ^{2}dt<\infty ; \\
\text{(iv)} & \left| f(t,y,z)-f(t,y,z^{\prime })\right| \leq C\left|
z-z^{\prime }\right| \text{, }\forall (t,y)\in [0,T]\times
\mathbb{R}\text{,
}z,z^{\prime }\in \mathbb{R}^{d}\text{, a.s.} \\
\text{(v)} & (y-y^{\prime })(f(t,y,z)-f(t,y^{\prime },z))\leq \mu
(y-y^{\prime })^{2}\text{, }\forall (t,z)\in [0,T]\times \mathbb{R}^{d}\text{%
, }y,y^{\prime }\in \mathbb{R}\text{, a.s.} \\
\text{(vi)} & y\rightarrow f(t,y,z)\text{ is continuous, }\forall
(t,z)\in [0,T]\times \mathbb{R}^{d}\text{, a.s.}
\end{array}
$
\\[0.1cm]
\textbf{(H3)} a barrier $(L_{t})_{0\leq t\leq T}$, which is a
continuous progressively measurable real-valued process, satisfying

\[
E[\varphi ^{2}(\sup_{0\leq t\leq T}(e^{\mu t}L_{t}^{+}))]<\infty ,
\]
and $(L_{t}^{+})_{0\leq t\leq T}\in \mathbf{S}^{2}(0,T)$, $L_{T}\leq
\xi $, a.s.

We shall  give an a priori estimate of the solution $(Y,Z,K)$ with
respect to the terminal condition $\xi $, the coefficient $f$ and
the barrier $L$. Unlike the Lipshitz case, we have in addition the
term $ E\varphi ^2(\sup_{0\leq t\leq T}(L_t^{+}))$ and a constant,
which only depends on $\varphi $, $\mu $, $k$ and $T$ :
\begin{proposition}
\label{esti}There exists a constant $C$, which only depends on
$T$, $\mu $ and $k$, such that
\begin{equation*}
\begin{split}
E \, \big[ \sup_{0\leq t\leq T}\left| Y_t\right| ^2 +\int_0^T\left|
Z_s\right| ^2ds+\left| K_T\right| ^2 \, \big]
&\leq C E \, \big [\xi ^2+\int_0^T f^2(t,0,0)dt + \varphi ^2(\sup_{0\leq t\leq T}(L_t^{+})) \big] \\
& \quad \quad +  C E [ \sup_{0\leq
t\leq T}(L_t^{+})^2+1+\varphi ^2(2T)].\\
\end{split}
\end{equation*}
\end{proposition}
\textbf{Proof.} Applying It\^{o}'s formula to $\left| Y_{t}\right|
^{2}$, and taking expectation, then
\begin{eqnarray*}
E[\left| Y_{t}\right| ^{2}+\int_{t}^{T}\left| Z_{s}\right| ^{2}ds]
&=&E[\left| \xi \right|
^{2}+2\int_{t}^{T}Y_{s}f(s,Y_{s},Z_{s})ds+2\int_{t}^{T}L_{s}dK_{s} \\
&\leq &E[\left| \xi \right|
^{2}+2\int_{t}^{T}Y_{s}f(s,0,0)ds+2\int_{t}^{T}(\mu \left|
Y_{s}\right| ^{2}+k\left| Y_{s}\right| \left| Z_{s}\right|
)ds+2\int_{t}^{T}L_{s}dK_{s}].
\end{eqnarray*}
It follows that
\begin{eqnarray*}
E[\left| Y_{t}\right| ^{2}+\frac{1}{2}\int_{t}^{T}\left|
Z_{s}\right| ^{2}ds] \leq E[\left| \xi \right|
^{2}+2\int_{t}^{T}f^{2}(s,0,0)ds+(2\mu
+1+2k^{2})\int_{t}^{T}\left| Y_{s}\right|
^{2}ds+2\int_{t}^{T}L_{s}dK_{s}].
\end{eqnarray*}
Then by Gronwall's inequality, we have
\begin{equation}
E\left| Y_{t}\right| ^{2}\leq CE[\left| \xi \right|
^{2}+\int_{0}^{T}f^{2}(s,0,0)ds+\int_{0}^{T}L_{s}dK_{s}],
\label{estu-y}
\end{equation}
then
\begin{equation}
E\int_{0}^{T}\left| Z_{s}\right| ^{2}ds\leq CE[\left| \xi \right|
^{2}+\int_{0}^{T}f^{2}(s,0,0)ds+\int_{0}^{T}L_{s}dK_{s}],
\label{estu-z}
\end{equation}
where $C$ is a constant only depends on $\mu $, $k$ and $T$, in
the following this constant can be changed line by line.

Now we estimate $K$ by approximation. By the existence of the
solution, theorem 2.2 in Lepeltier et al. \cite{LMX}, we take the
process $Z$ as a known process. Without losing generality we write
$f(t,y)$ for $f(t,y,Z_t)$, here $f(t,0)=f(t,0,Z_t) $ is a process in
$\mathbf{H}^2(0,T)$. Set
\begin{eqnarray*}
\xi ^{m,n} &=&(\xi \vee (-n))\wedge m, \\
f^{m,n}(t,y) &=&f(t,y)-f(t,0)+(f(t,0)\vee (-n))\wedge m.
\end{eqnarray*}
For $m$, $n\in \mathbf{N}$, $\xi ^{m,n}$ and $\sup_{0\leq t\leq
T}f^{m,n}(t,0)$ are uniformly bounded. Consider the RBSDE$(\xi
^{m,n},f^{m,n},L)$,
\begin{eqnarray*}
Y_t^{m,n} &=&\xi
^{m,n}+\int_t^Tf^{m,n}(t,Y_s^{m,n})ds+K_T^{m,n}-K_t^{m,n}-%
\int_t^TZ_s^{m,n}dB_s, \\
Y_t^{m,n} &\geq &L_t,\int_0^T(Y_s^{m,n}-L_s)dK_s^{m,n}=0.
\end{eqnarray*}
if we recall the transform in step 2 of the proof of theorem 2.2 in
Lepeltier et al. \cite{LMX},
since $\xi ^{m,n}$, $f^{m,n}(t,0)\leq m$, we know that $%
(Y_t^{m,n},Z_t^{m,n},K_t^{m,n})$ is the solution of this RBSDE, if
and only
if $(Y^{m,n\prime },Z^{m,n\prime },K^{m,n\prime })$ is the solution of RBSDE$%
(\xi ^{m,n\prime },f^{m,n\prime },L^{\prime })$, where
\begin{equation*}
(Y_t^{m,n\prime }, Z_t^{m,n\prime }, K_t^{m,n\prime })=
(Y_t^{m,n}+m(t-2(T\vee 1)), Z_t^{m,n}, K_t^{m,n})
\end{equation*}
and
\begin{eqnarray*}
\xi ^{m,n\prime } &=&\xi ^{m,n}+2mT-m(T\vee 1), \\
f^{m,n\prime }(t,y) &=&f^{m,n}(t,y-m(t-2(T\vee 1)))-m, \\
L_t^{\prime } &=&L_t+m(t-2(T\vee 1)).
\end{eqnarray*}
Without losing generality we set $T\geq 1$. Then $\xi ^{m,n\prime
}\leq 0$ and $f^{m,n\prime }(t,0)\leq 0$. Since $(Y^{m,n\prime
},Z^{m,n\prime },K^{m,n\prime })$ is the solution of RBSDE$(\xi
^{m,n\prime },f^{m,n\prime },L^{\prime })$, then we have
\[
K_T^{m,n\prime }=Y_0^{m,n\prime }-\xi ^{m,n\prime
}-\int_0^Tf^{m,n\prime }(s,Y_s^{m,n\prime
},Z_s)ds+\int_0^TZ_s^{m,n\prime }dB_s,
\]
which follows
\begin{equation}
\label{est-k-bon} E[(K_T^{m,n\prime })^2]\leq 4E[\left|
Y_0^{m,n\prime }\right| ^2+\left| \xi ^{m,n\prime }\right|
^2+(\int_0^Tf^{m,n\prime }(s,Y_s^{m,n\prime })ds)^2+\int_0^T\left|
Z_s^{m,n\prime }\right| ^2ds].
\end{equation}
Applying It\^o's formula to $\left| Y^{m,n}\right| ^2$, like
(\ref{estu-y}) and (\ref{estu-z}), we have
\[
E\left| Y_t^{m,n}\right| ^2+E\int_t^T\left| Z_s^{m,n}\right|
^2ds\leq CE[\left| \xi ^{m,n}\right|
^2+\int_t^T(f^{m,n}(s,0))^2ds+\int_t^TL_sdK_s^{m,n}].
\]
So
\begin{eqnarray*}
\left| Y_0^{m,n\prime }\right| ^2+\int_0^T\left| Z_s^{m,n\prime
}\right| ^2ds &=&2\left| Y_0^{m,n}\right|
^2+8m^2T^2+E\int_0^T\left| Z_s^{m,n}\right|
^2ds \\
&\leq &CE[\left| \xi ^{m,n}\right|
^2+\int_0^T(f^{m,n}(s,0)ds)^2+\int_0^TL_sdK_s^{m,n}]+8m^2T^2.
\end{eqnarray*}
For the third term on the right side of (\ref{est-k-bon}), from
Lemma 2.3 in Lepeltier et al. \cite{LMX}, we remember that
\begin{equation}
(\int_0^Tf^{m,n\prime }(s,Y_s^{m,n\prime })ds)^2\leq \max
\{(\int_0^Tf^{m,n\prime }(s,\widetilde{Y}_s^{m,n})ds)^2,(\int_0^Tf^{m,n%
\prime }(s,\overline{Y}_s^{m,n})ds)^2\},  \label{control-f}
\end{equation}
where $(\widetilde{Y}^{m,n},\widetilde{Z}^{m,n})$ is the solution
the following BSDE
\begin{equation}
\widetilde{Y}_t^{m,n}=\xi ^{m,n\prime }+\int_0^Tf^{m,n\prime }(s,\widetilde{Y%
}_s^{m,n})ds-\int_0^T\widetilde{Z}_s^{m,n}dB_s,  \label{con-BSDE}
\end{equation}
and
\[
\overline{Y}_s^{m,n}=ess\sup_{\tau \in \mathcal{T}_{t,T}}E[(L_\tau
^{^{\prime }})^{+}1_{\{\tau <T\}}+(\xi ^{m,n})^{+}1_{\{\tau =T\}}|\mathcal{F}%
_t]=ess\sup_{\tau \in \mathcal{T}_{t,T}}E[(L_\tau ^{^{\prime }})^{+}|%
\mathcal{F}_t].
\]
From (\ref{con-BSDE}), and proposition 2.2 in Pardoux
\cite{Pardoux99}, we have
\begin{eqnarray*}
E(\int_0^Tf^{m,n\prime }(s,\widetilde{Y}_s^{m,n})ds)^2 &\leq
&CE[\left| \xi
^{m,n\prime }\right| ^2+(\int_0^Tf^{m,n\prime }(s,0)ds)^2] \\
&\leq &CE[\left| \xi ^{m,n}\right|
^2+\int_0^T(f^{m,n}(s,0))^2ds]+C\varphi ^2(2mT)+Cm^2.
\end{eqnarray*}
For the other term in (\ref{control-f}), with $\sup_{0\leq t\leq
T}\overline{Y}_s^{m,n}=\sup_{0\leq t\leq T}(L_t^{\prime })^{+}$,
we get
\begin{eqnarray*}
E(\int_0^Tf^{m,n\prime }(s,\overline{Y}_s^{m,n})ds)^2 &\leq
&E[\int_0^T2(f^{m,n\prime }(s,0))^2ds+2T\varphi ^2(\sup_{0\leq
t\leq
T}(L_t^{\prime })^{+}] \\
&\leq &E[4\int_0^Tf^{m,n}(s,0)^2ds+2T\varphi ^2(\sup_{0\leq t\leq
T}(L_t)^{+})]+2m^2T+4T\varphi ^2(2mT).
\end{eqnarray*}
Consequently, we deduce that
\begin{eqnarray*}
E[(K_T^{m,n})^2] &=&E[(K_T^{m,n\prime })^2] \\
&\leq &CE[\left| \xi ^{m,n}\right|
^2+\int_0^T(f^{m,n}(s,0))^2ds+\int_t^TL_sdK_s^{m,n}+\varphi
^2(\sup_{0\leq
t\leq T}(L_t)^{+})+m^2+\varphi ^2(2mT)] \\
&\leq &CE[\left| \xi \right| ^2+\int_0^T(f(s,0,Z_s))^2ds+\varphi
^2(\sup_{0\leq t\leq T}(L_t)^{+})+\sup_{0\leq t\leq
T}((L_t)^{+})^2]+\frac
12E[(K_T^{m,n})^2] \\
&&+C(m^2+\varphi ^2(2mT)).
\end{eqnarray*}
Moreover using (\ref{estu-z}) and the fact that $f$ is Lipschitz on
$z$, it follows that
\begin{eqnarray}
E[(K_T^{m,n})^2] &\leq &CE[\left| \xi \right|
^2+\int_0^T(f(s,0,0))^2ds+\varphi ^2(\sup_{0\leq t\leq
T}(L_t)^{+})+\sup_{0\leq t\leq T}((L_t)^{+})^2 \label{est-k-mn}\\
&&+\int_0^TL_sdK_s]+C(m^2+\varphi ^2(2mT)). \nonumber
\end{eqnarray}

Let $m\rightarrow \infty $, then
\[
E[\left| \xi ^{m,n}-\xi ^{n}\right| ^{2}]\rightarrow
0,E\int_{0}^{T}\left| f^{m,n}(t,0)-f^{n}(t,0)\right|
^{2}\rightarrow 0,
\]
where $\xi ^{n} =\xi \vee (-n)$ and $f^{n}(t,y)
=f(t,y)-f(t,0)+f(t,0)\vee (-n)$.

Thanks to the convergence result of step 3 of the proof for theorem 2.2 in \cite{LMX}, we know that $%
(Y^{m,n},Z^{m,n},K^{m,n})\rightarrow (Y^{n},Z^{n},K^{n})$ in $\mathbf{S}%
^{2}(0,T)\times \mathbf{H}_{d}^{2}(0,T)\times \mathbf{A}^{2}(0,T)$, where $%
(Y^{n},Z^{n},K^{n})$ is the soultion of the RBSDE$(\xi
^{n},f^{n},L)$.
Moreover $K_{T}^{m,n}\searrow K_{T}^{n}$ in $\mathbf{L}^{2}(\mathcal{F}_{T})$%
, so we have $K_{T}^{n}\leq K_{T}^{1,n}$, which implies for each
$n\in \mathbf{N}$,
\begin{eqnarray}
E[(K_{T}^{n})^{2}] \leq E[(K_{T}^{1,n})^{2}]\label{est-k-n}
\end{eqnarray}

Then, letting $n\rightarrow \infty $, by the convergence result in
step 4, since
\[
E[\left| \xi ^n-\xi \right| ^2]\rightarrow 0,E\int_0^T\left|
f^n(t,0)-f(t,0)\right| ^2\rightarrow 0,
\]
the sequence $(Y^n,Z^n,K^n)\rightarrow (Y,Z,K)$ in
$\mathbf{S}^2(0,T)\times \mathbf{H}_d^2(0,T)\times
\mathbf{A}^2(0,T)$, where $(Y,Z,K)$ is the solution of the
RBSDE$(\xi ,f,L)$. From (\ref{est-k-n}), and (\ref{est-k-mn}) for
$m=1$, we get
\begin{eqnarray*}
\ E[(K_T)^2]\  &\leq &CE[\left| \xi \right|
^2+\int_0^T(f(s,0,0))^2ds+ \varphi ^2(\sup_{0\leq t\leq
T}(L_t)^{+}) + \sup_{0\leq t\leq
T}((L_t)^{+})^2 \\
&& +\int_0^TL_sdK_s]+C(1+\varphi ^2(2T))\\
&\leq &CE[\left| \xi \right| ^2+\int_0^T(f(s,0,0))^2ds+\varphi
^2(\sup_{0\leq t\leq T}(L_t)^{+})+\sup_{0\leq t\leq T}((L_t)^{+})^2] \\
&&\ +\frac 12E[(K_T)^2]+C(1+\varphi ^2(2T)).
\end{eqnarray*}
Then it follows that for each $t\in [0,T]$,
\begin{eqnarray*}
E[\left| Y_t\right| ^2+\int_0^T\left| Z_s\right| ^2ds+(K_T)^2]
&\leq &CE[\left| \xi \right| ^2+\int_0^T(f(s,0,0))^2ds+\varphi
^2(\sup_{0\leq t\leq T}(L_t)^{+}) \\
&&\ +\sup_{0\leq t\leq T}((L_t)^{+})^2]+C(1+\varphi ^2(2T)).
\end{eqnarray*}
Finally we get the result, by applying BDG inequality. $\hfill \Box$

\end{document}